\documentclass[12pt]{article} 
\usepackage{amsmath,amssymb,eucal,bbm,mathtools} 
\font\tenrm=cmr10

\font\cmssl=cmss10 at 12 pt  
   
\font\bigss=cmssdc10 scaled 2300

\font\cmsslll=cmss10 at 14 pt

\renewcommand{\a}{\alpha}  
\renewcommand{\b}{\beta}  
  
\renewcommand{\d}{\delta}  
\newcommand{\e}{\epsilon}  
  
\newcommand{\g}{\gamma}

\renewcommand{\k}{\kappa}  
\renewcommand{\l}{\lambda}  
\renewcommand{\o}{\omega}

\newcommand{\x}{\xi}

\newcommand{\G}{\Gamma}

\newcommand{\bC}{\mathbb{C}}  
\newcommand{\bR}{\mathbb{R}}  
  
\newcommand{\bH}{\mathbb{H}}

\newcommand{\ga}{\mathfrak{a}}

\newcommand{\gl}{\mathfrak{l}}  
  
\newcommand{\gn}{\mathfrak{n}}

\newcommand\Sp{\mathrm{Sp}}

\newcommand\SO{\mathrm{SO}}  
\newcommand\SU{\mathrm{SU}}  
\newcommand\U{\mathrm{U}}  
  
\renewcommand{\square}{\kern1pt\vbox  
             {\hrule height 0.6pt\hbox{\vrule width 0.6pt\hskip 3pt  
  \vbox{\vskip 6pt}\hskip 3pt\vrule width 0.6pt}\hrule height0.6pt}  
  \kern1pt}  

\newcommand{\ra}{\rightarrow}

\DeclareMathOperator\End{End\;}

\newtheorem{Th}{Theorem}  
\newtheorem{Prop}{Proposition}  
  
\newtheorem{Cor}{Corollary}  
\newtheorem{Lem}{Lemma}  
\newtheorem{Def}{Definition}  
\newcommand{\bt}{\begin{Th}\ \ }  
\newcommand{\et}{\end{Th}}  
\newcommand{\bp}{\begin{Prop}\ \ }  
\newcommand{\ep}{\end{Prop}}  
\newcommand{\bc}{\begin{Cor}\ \ }  
\newcommand{\ec}{\end{Cor}}  
\newcommand{\bl}{\begin{Lem}\ \ }  
\newcommand{\el}{\end{Lem}}  
\newcommand{\bd}{\begin{Def}\ \ }  
\newcommand{\ed}{\end{Def}}  

\newcommand{\pf}{\noindent{\it Proof:\ \ }}  
\newcommand{\qed}{\hfill\square}  
\newcommand{\n}{\nabla}  
  
\newcommand{\ot}{\otimes}

\newcommand{\be}{\begin{equation}}  
\newcommand{\ee}{\end{equation}}  
  
\newcommand\re[1]{(\ref{#1})}  
\newcommand{\arr}{\begin{array}{rlll}}  
\newcommand{\ea}{\end{array}}  
\newcommand{\bea}{\begin{eqnarray}}  
\newcommand{\eea}{\end{eqnarray}}  
\newcommand{\bean}{\begin{eqnarray*}}  
\newcommand{\eean}{\end{eqnarray*}}  
  
\catcode`@=11  
\@addtoreset{equation}{section}  
\catcode`@=12

\newcommand\iu{\operatorname{i}}
\newcommand\diff{\mathrm{d}}
\renewcommand{\Re}{\operatorname{Re}}
\renewcommand{\Im}{\operatorname{Im}}

\def\cN{{\cal N}}

\begin{document}  
\rightline{IPhT-t11/196} 
\rightline{\small ZMP-HH/11-18}
\vskip 1.5 true cm  
\begin{center}  
{\bigss On certain K\"ahler quotients\\[1em] of quaternionic K\"ahler manifolds}\\[.5em]
\vskip 1.0 true cm   
{\cmsslll  V.\ Cort\'es$^{1,2}$, J.\ Louis$^{2,3}$, P.\ Smyth$^{4}$ 
and H.\ Triendl$^{5}$} \\[5pt] 
$^1$ {\tenrm   Department Mathematik, Universit\"at Hamburg\\ 
Bundesstra{\ss}e 55, 
D-20146 Hamburg, Germany\\  
cortes@math.uni-hamburg.de}\\[1em]   
$^2$ {\tenrm   Zentrum f\"ur Mathematische Physik, Universit\"at Hamburg\\ 
Bundesstra{\ss}e 55, 
D-20146 Hamburg, Germany}\\[1em]   
$^3$ {\tenrm II.\ Institut f\"ur Theoretische Physik, Universit\"at Hamburg\\ 
Luruper Chaussee 149, D-22761 Hamburg, Germany\\ 
jan.louis@desy.de}\\[1em]
$^4$ {\tenrm Institut de Th\'eorie des Ph\'enom\`enes Physiques, EPFL\\ CH-1015 Lausanne, Switzerland\\paul.smyth@epfl.ch}\\[1em]   
$^5$ {\tenrm Institut de Physique Th\'eorique, CEA Saclay\\ 
Orme des Merisiers, F-91191 Gif-sur-Yvette, France\\ 
hagen.triendl@cea.fr} 
\end{center}  
\vskip 1.0 true cm  
\baselineskip=18pt  
\begin{abstract}  
\noindent  
We prove that, given a certain isometric action of a two-dimensional 
Abelian group $A$ on a quaternionic K\"ahler manifold $M$ which preserves
a submanifold $N\subset M$, the quotient $M'=N/A$ has a natural 
K\"ahler structure. 
We verify that the assumptions on the group action 
and on the submanifold $N\subset M$ are 
satisfied for a large class of examples obtained 
from the supergravity c-map.  In particular, we find that all quaternionic
K\"ahler manifolds $M$ in the image of the c-map admit an integrable
complex structure compatible with  the quaternionic structure, such that 
$N\subset M$ is a complex submanifold.  
Finally, we discuss how
the existence of the K\"ahler structure on $M'$ is required 
by the consistency of spontaneous $\cN=2$ to $\cN=1$ 
supersymmetry breaking.
\end{abstract}

\section*{Introduction} 
 
Since the work of Galicki and Lawson \cite{GL} it has been known that 
a quaternionic analogue of the well-known symplectic reduction exists. 
In fact, as shown in \cite{ACDV}, any isometric 
action of a Lie group $G$ on 
a quaternionic K\"ahler manifold $(M,g,Q)$  of nonzero scalar 
curvature gives rise to a $\mathfrak{g}^*$-valued \
section $P\in \G (Q \ot \mathfrak{g}^*)$ of the 
quaternionic structure $Q\subset \End (TM)$. 
$P$ is called the moment map and by taking the
quotient $\{ P= 0\}/G$ one obtains a new quaternionic K\"ahler 
manifold, provided that the usual regularity assumptions are fulfilled.  

In this paper, however,  
we are interested in constructing K\"ahler manifolds
out of  quaternionic K\"ahler manifolds. Such a procedure is needed in order   
to break supersymmetry from $\cN=2$ to $\cN=1$ in supersymmetric theories of 
gravity in four spacetime dimensions  \cite{Ferrara:1995gu,JL,LST1,LST2}. The reason is that 
quaternionic K\"ahler manifolds of negative scalar curvature 
occur as scalar manifolds of $\cN=2$ supergravity, whereas 
$\cN=1$ supergravity requires the 
scalar manifold to be K\"ahler. A natural but rather restrictive way to relate
quaternionic K\"ahler manifolds to K\"ahler manifolds of lower dimension
is to consider \emph{K\"ahler submanifolds} 
$(N,g_N,J_N) \subset (M,g,Q)$, 
such that $g_N = g|_N$ and $J_N\in \G (Q|_N)$. It is shown in \cite{AM}
that the dimension of such a submanifold cannot exceed $2n$ if 
$M$ has nonzero scalar curvature, where $\dim M = 4n$. 

Our new idea is to drop the K\"ahler condition on $(N,g_N,J_N)$ 
still maintaining the integrability of $J_N\in \G (Q|_N)$. 
The final K\"ahler manifold $M'$ is then obtained as an 
appropriate quotient $M'=N/A$ of $N$. To define the quotient
we make use of two commuting Killing vector fields $\xi_1$, $\xi_2$,
which generate a free proper isometric action
of a two-dimensional Lie group $A$.  The necessary technical assumptions
on $\xi_1, \xi_2$ for our construction 
are formulated in terms of the corresponding moment maps
$P_1, P_2\in \G(Q)$, see Theorem \ref{thm:main} and Corollary \ref{cor:main}.  
The main result can be summarized as follows.
\bt Let $M$ be a quaternionic K\"ahler manifold of nonzero scalar curvature,
$N\subset M$ a submanifold and $\xi_1, \xi_2$ Killing vector fields of $M$ 
which satisfy the assumptions of Theorem \ref{thm:main} and Corollary 
\ref{cor:main}. Then $M'=N/A$ carries an induced K\"ahler structure,
where $A$ is the transformation group generated by $\xi_1, \xi_2$. 
\et 

The main body of the article is devoted to the investigation
of several classes of examples. As a first and simplest example we take 
$N=M=H_\bR^4=H_\bH^1$ the real hyperbolic four-space (which coincides 
with quaternionic hyperbolic line) and obtain $M'=H_\bC^1$. Then
we study the quaternionic K\"ahler manifolds $(M,g,Q)$ in the image of
the c-map \cite{Cecotti:1988qn,Ferrara:1989ik}. These manifolds 
have negative scalar curvature and 
are associated with a (projective) special K\"ahler domain $M_{\rm sk}$, 
the geometry of which can be encoded in a holomorphic prepotential $F(Z)=F(Z^1,\ldots ,Z^n)$.  As a first step in the study of the c-map examples we 
obtain the following general result, see Proposition \ref{J3Prop} and 
\ref{biholProp}. 
\bt Let $(M,g,Q)$ be a quaternionic K\"ahler manifold in the image
of the c-map. Then the quaternionic structure $Q$ of $M$ admits 
a global orthonormal frame $(J_1,J_2,J_3)$ such that the almost
complex structure $J_3\in \G (Q)$ is integrable. $(M,J_3)$ is the total
space of a holomorphic submersion $M\ra M_{\rm sk}$ with all fibers
biholomorphic to the domain $\bR^{n+1} + i V \subset \bC^{n+1}$, where
$\dim M = 4n$ and 
\be \label{indEqu} V= \{ (x_0,x_1,\ldots ,x_n)\in \bR^{n+1}| x_0 > 
\sum_{i=1}^{n-1} x_i^2 -x_n^2 \}.\ee
\et 
This should be contrasted 
with the situation for complete quaternionic K\"ahler manifolds of positive
scalar curvature, which do not even admit an \emph{almost} 
complex structure compatible\footnote{Note that the complex 
Grassmannians $Gr_2(\bC^n)$ ($n\ge 3$) do admit a  
complex structure which is even K\"ahler for the quaternionic K\"ahler metric 
but it does not belong to the quaternionic structure. It is known
that these complex Grassmannians are the only complete quaternionic K\"ahler manifolds 
of positive scalar curvature which admit an almost complex structure 
\cite{GMS}.}   
with the quaternionic structure \cite{AMP}. 
Some interesting properties of the complex structure $J_3$ are described 
in Proposition \ref{J3Prop} and Proposition 
\ref{biholProp}. 
We then define a complex submanifold 
\[ N \subset (M,J_3),\] 
see 
Proposition \ref{NProp},  
associated with a choice of a null vector $v_0 \in TM_{\rm ask}$, 
where $M_{\rm ask} \ra M_{\rm sk}$ is the affine
special K\"ahler manifold associated with $M_{\rm sk}$. 
The complex codimension of $N\subset M$ is $r+1$, where 
$r$ is the rank of a certain complex matrix $(G_{AB})$,  which depends on 
the choice of $v_0$, see equation\footnote{Note that 
$v_0=\sum D^A \partial/\partial Z^A|_{Z_0} +$ c.c. .} \re{equ:GAB} and the remark 
on page \pageref{GAB}. 
The structure of the complex manifold $N$ is described in Proposition
\ref{NbiholProp}. In particular, we find that $N$ is always the total
space of a holomorphic submersion 
\[ N\ra  M_{\rm sk}^\wedge,\] 
where
$M_{\rm sk}^\wedge\subset M_{\rm sk}$ is a complex submanifold 
and the fibers are biholomorphic to $B^{n-1}_\bC\times \bC$. The inclusion
$N\subset M$ maps the fibers of $N\ra M_{\rm sk}^\wedge$ into the fibers of
$M\ra M_{\rm sk}$. 
Next we define two Killing vector fields
$\xi_1, \xi_2$ on $M$, which depend on the choice of $v_0$. 
We show in Proposition \ref{holactionProp}  that they are tangent to 
$N\subset M$ and generate a holomorphic, free and proper 
action of the additive group $A=\bC$ on $N$.  
We then have the following result, cf.\ Theorem \ref{ThmM'}. 
\bt The resulting quotient $M'=N/A$ is always  the total space of 
a holomorphic submersion
\[ M'\ra M^{\wedge}_{\rm sk},\]
where the fibers are isomorphic to the complex ball $B^{n-1}_{\bC}\cong 
\bC H^{n-1}$ with its standard complex hyperbolic metric 
of constant holomorphic sectional
curvature $-4$. The projection $N\ra M'= N/A$ maps the fibers of
$N\ra M^{\wedge}_{\rm sk}$ to the fibers of $M'\ra M^{\wedge}_{\rm sk}$. 
\et 
We also show that $M'$ is complete if and only if the base manifold 
$M_{\rm sk}^\wedge$ is complete, see Remark on page 
\pageref{completenesspageref}.  
Let us emphasize a subtle but crucial point in the construction. 
The fibers $M_p=\pi^{-1}(p)$ of $\pi : (M,g,J_3)\ra M_{\rm sk}$ consist of a 
solvable Lie group $G$ endowed with a family of left-invariant metrics 
$g_G(p)$ and left-invariant skew-symmetric 
complex structures $J_G(p)$:
\[ (M_p,g|_{M_p},J_3|_{M_p}) = (G,g_G(p),J_G(p)),\quad p\in M_{\rm sk}.\] 
The group $G$ is precisely the 
Iwasawa subgroup of $\SU(1,n+1)$, which is the group
of holomorphic isometries of the complex hyperbolic space  
$\bC H^{n+1}=\U(1,n+1)/(\U(1)\times \U(n+1))=\SU(1,n+1)/
\mathrm{S}(\U(1)\times \U(n+1))$. 
Since $G$ acts simply transitively on $\bC H^{n+1}$, 
we can identify the K\"ahler manifold $\bC H^{n+1}$ with $(G,g_{can},J_{can})$,
where $(g_{can},J_{can})$ is a left-invariant K\"ahler structure on $G$:
\[ \bC H^{n+1} = (G,g_{can},J_{can}).\]  
{}From the Riemannian point of view, the fibers   
$(M_p,g|_{M_p})=(G,g_G(p))$ are as nice as possible. They are all 
isometric to  $(G,g_{can})\cong \bC H^{n+1}$, although the metric $g_G(p)$ 
is never independent of $p\in M_{\rm sk}$. 
However,  in view of the above discussion, the Hermitian manifold 
$(G,g_G(p),J_G(p))$ cannot be K\"ahler, 
since $2n+2=\dim G > \frac{1}{2}\dim M=2n$. This means that
$J_G(p)$ does not coincide with the canonical (parallel) complex 
structure $J_{can}(p)$ on $(G,g_G(p))H^{n+1}$, for which 
$(G,g_G(p),J_{can}(p))\cong \bC H^{n+1}=(G,g_{can},J_{can})$. 
One can show that $(G,J_G(p))$ is not even
biholomorphic\footnote{A proof of this fact can be found in 
\cite{CH}, which includes the classification
of skew-symmetric left-invariant complex structures on 
$(G,g_{can})=\bC H^{n+1}$.} to $\bC H^{n+1}$. This is related to 
the non-positivity of the quadratic form  on the right-hand side of the 
inequality defining the complex domain $\bR^{n+1} + i V \subset \bC^{n+1}
\cong (G,J_G(p))$, see \re{indEqu}.  Summarizing, we have that
\[ (G,g_G(p))\cong (G,g_{can})\cong \bC H^{n+1}\quad \mbox{but}\quad 
(G,J_G(p))\not\cong (G,J_{can})\cong \bC H^{n+1}.\]  
It turns out that when passing to the quotient $M'$, the fibers 
$M_p'$, $p\in M_{\rm sk}^\wedge$, 
become all isometric \emph{and} biholomorphic to 
$\bC H^{n-1}$. In fact, we show that by considering the submanifold 
$N_p := N\cap M_p\subset M_p$, which is biholomorphic to 
$\bC H^{n-1}\times \bC$, and its quotient $M_p'=N_p/A$ we 
reduce the domain $M_p\cong \bR^{n+1}+iV\not\cong \bC H^{n+1}$ to 
$M_p'\cong (\bR^{n+1}+iV)\cap \bC^{n-1}=\bR^{n-1}+iV'$, where now
\[ V'  = \{ (x_0,x_1,\ldots ,x_{n-2})\in \bR^{n-1}| x_0 > 
\sum_{i=1}^{n-2} x_i^2  \}\]   
is defined by a positive definite quadratic form. Therefore, 
$\bR^{n-1}+iV'$ is biholomorphic to $\bC H^{n-1}$. 

More detailed information is obtained in 
Section \ref{QuadprepotSec} and Section \ref{sect:cubic} when
the prepotential is either quadratic or of the form 
$F=\frac{h(Z^1,\ldots ,Z^{n-1})}{Z^n}$, where $h$ is a homogeneous cubic
polynomial with real coefficients. As usual in the physics literature, 
the latter class will be simply referred to as having \emph{cubic 
prepotential}.  It is particularly 
interesting for string theory compactifications and contains a wealth of
homogeneous as well as inhomogeneous examples. We show that in the case
of cubic prepotential the dimension of $M'$ can be as large as $\dim M -8$
with $\dim M$ arbitrarily big.  
In the case of quadratic prepotential the structure of 
$M'$ is completely determined as follows, cf.\ 
Theorem \ref{quadrThm}.
\bt The  K\"ahler manifolds $M'$ obtained from the above 
quotient construction  applied to 
the quaternionic K\"ahler manifold $M = \frac{\U(2,n)}{\U(2)\times \U(n)}\ra 
M_{\rm sk}=H_\bC^{n-1}$  are always isomorphic to $H_\bC^{n-1}\times H_\bC^{n-1}$, 
provided that $M_{\rm sk}^\wedge\subset M_{\rm sk}$ is complete. 
In this case, the holomorphic submersion $M'\ra M_{\rm sk}^\wedge$ is trivial 
and $M_{\rm sk}^\wedge=M_{\rm sk}$.   
\et   
So in this case, $\dim M' = \dim M -4$. 

The mathematical results obtained in this paper
are motivated by the consistency of spontaneous $\cN=2$ to $\cN=1$
supersymmetry breaking \cite{Ferrara:1995gu,JL,LST1,LST2} and 
in Section~\ref{sect:physics} we briefly discuss this relation. 
Quaternionic K\"ahler manifolds appear
naturally in $\cN=2$ supergravity theories as part of the scalar field
space. The Higgs mechanism responsible for the supersymmetry breaking
requires two massive vector fields coupled to two 
Killing vector fields that fulfill the assumptions of Theorem \ref{thm:main}.
Furthermore,
an $\cN=1$ effective action can be defined below the scale of
supersymmetry breaking and is obtained by integrating out all 
massive degrees of freedom. 
Integrating out massive scalars 
corresponds to taking a
submanifold $N\subset M$, while integrating out two massive 
vector fields corresponds to taking the quotient with respect 
to the two-dimensional Abelian Lie group $A$ generated by the two
Killing vectors, as specified in Theorem~1.
 Consistency with $\cN=1$ supersymmetry implies that the resulting scalar field space  $M'=N/A$ should be K\"ahler. 
\section{Basic results about quaternionic K\"ahler manifolds}
In this section we recall some known facts about quaternionic K\"ahler
manifolds, see e.g.\ \cite{ACDV} for more details. 
\bd A {\cmssl quaternionic K\"ahler manifold} $(M,g,Q)$ 
is a Riemannian manifold $(M,g)$ which is endowed with a parallel
skew-symmetric  quaternionic structure $Q\subset \End TM$. 
If $\dim M =4$ we require, in addition, that $Q\cdot R=0$. 
(This condition is automatically satisfied if $\dim M >4$.)
\ed
Let $(J_\a)_{\a =1,2,3}$ be an orthonormal local frame of $Q$
such that $J_3=J_1J_2$. Then 
\be \label{e:omega} \n J_\a = -(\o_\b \ot J_\g -\o_\g \ot J_\b ),\ee
for some triplet of connection forms $\o_\a$, where $(\a ,\b ,\g )$ 
is always a cyclic permutation of $(1,2,3)$. These one-forms are related
to the fundamental two-forms $\varphi_\a = g(\cdot , J_\a)$ by the 
following structure equations: 
\be \label{e:str}\nu \varphi_\a = d\o_\a + \o_\b\wedge \o_\g ,\ee
where $\nu = \frac{scal}{4n(n+2)}$ 
stands for the reduced scalar curvature, the quotient of the scalar 
curvature of $(M,g)$ by that of $\bH P^n$, with $4n=\dim M$. 
Quaternionic K\"ahler manifolds are Einstein; in particular, 
$\nu$ is a constant. 

Now let $\xi$ be a Killing vector field on a quaternionic K\"ahler manifold
of nonzero scalar curvature, i.e.\ $\nu\neq 0$. Then $Q$ is invariant under the
flow of $\xi$, as well as under parallel transport. This implies that
the endomorphism field $\n \xi$ is a section of the normaliser
\[ N(Q)=Q\oplus Z(Q)\] 
of $Q$ in $\mathfrak{so}(TM)$.  Here $Z(Q)$ stands for the centraliser
of $Q$. Note that 
\[ N(Q)_p\cong \mathfrak{sp}(1)\oplus \mathfrak{sp}(n)\quad \forall p\in M,\]
where $\mathfrak{sp}(n)$ is the Lie algebra of the compact
symplectic group $\Sp (n)$, which is usually denoted by $\mathrm{USp}(2n)$ 
in the physics literature. Let us use 
\be \label{momentEqu} P := (\n \xi )^Q\in \G (Q)\ee 
to denote the projection of $\n \xi$ onto $Q$.  The section $P : M \ra Q$ is called
the {\cmssl moment map} associated with $\xi$. Its covariant
derivative is given by:     
\be \label{e:moment}
\n P = \frac{\nu}{2} \sum \varphi_\a (\cdot ,\xi )\ot J_\a .\ee
For the last formula, see \cite{ACDV} Proposition 2. 

\section{The new quotient construction} \label{sect:quotient}
\bt \label{thm:main} Let $(M, g, Q)$ be a quaternionic K\"ahler 
manifold of nonzero scalar curvature, $\xi_1$, $\xi_2$ two 
Killing vector fields with corresponding moment maps $P_i\in \G (Q)$, $i=1, 2$, $N\subset M$ a submanifold and $\mathfrak{X}(N)$ the space of smooth vector fields on $N$ such that:
\begin{enumerate}
\item[(i)] $\xi_1|_N$, $\xi_2|_N\in \mathfrak{X}(N)$,  
$[\xi_1 ,\xi_2]|_N =0$ and $|\xi_1|=|\xi_2|\neq 0$ on $N$, 
\item[(ii)] $P_1P_2|_N$ is a section of $Q|_N$ which preserves
$TN$ and maps $\xi_1|_N$ to $f\xi_2|_N$, where $f\in C^\infty (N)$ is some 
nowhere vanishing function. 
\end{enumerate}
Then the integrable distribution 
$\mathcal{D}\subset TN$ spanned by $\xi_1|_N$, $\xi_2|_N$
has an induced transversal K\"ahler structure $(h, J)$. 
The complex structure $J$ is induced by 
$I:=\frac{1}{f}P_1P_2|_N\in \G (N,Q)$, which defines an
integrable complex structure on $N$.   
\et  

\noindent
{\bf Remarks:} 1) 
We will show below that for the quaternionic K\"ahler manifolds
$(M,g,Q)$ in the image of the c-map there exists a global orthonormal
frame $(J_1,J_2,J_3)$ of $Q$ such that the almost complex structure $J_3$
is integrable. The above construction will then be applied to
an appropriate complex submanifold $N$ of $(M,J_3)$.\\
2) The quaternionic K\"ahler manifolds in the image of the c-map
include all the known homogeneous quaternionic K\"ahler manifolds 
of negative scalar curvature 
with the exception of the quaternionic hyperbolic spaces $H^n_\bH$, $n\ge 2$.

\pf 
Let us use $\mathcal{N} = \mathcal{D}^\perp\cong 
TN/\mathcal{D}$ to denote the Riemannian 
normal bundle of $\mathcal{D}$ in $N$. We then define
the transversal metric $h\in \G (S^2\mathcal{N}^*)$ 
as the restriction 
\[ h := g|_{\mathcal{N}\times \mathcal{N}}.\]
It follows from (i-ii) 
that $I:=\frac{1}{f}P_1P_2\in \G (N,Q)$ is an almost complex structure.
We can choose an orthonormal local frame $(J_\a)_{\a =1,2,3}$
of $Q$ such that $J_3=I$ on $N$ and $J_3=J_1J_2$. 
Since $I$ preserves 
$\mathcal{D}\subset TN$, $TN$ and therefore $\mathcal{N}=
\mathcal{D}^\perp\subset TN$, 
we can define 
\[ J := I|_{\mathcal{N}}\in \G (N,\End \mathcal{N}).\]
Clearly, $J_p$ is a skew-symmetric  complex structure 
on the Euclidian vector space $(\mathcal{N}_p,h_p)$, for all $p\in M$. 
We claim that $(h,J)$ defines a transversal K\"ahler structure for
the foliation of $N$ defined by the integral surfaces of the distribution 
$\mathcal{D}$. This means that $(h,J)$ induces a K\"ahler structure
on any submanifold $S\subset N$ transversal to $\mathcal{D}$ and that
the K\"ahler structures on a pair of such submanifolds 
$S, S'\subset M$ intersecting the same leaves are related
by the corresponding holonomy transformation of the foliation. 
To prove this it suffices to check that 
\bea \mathcal{L}_{\xi_i}g&=&0,\label{e:isom}\\
\mathcal{L}_{\xi_i}I&=&0,\label{e:I}\\
{[}X,Y{]} &\in&  \G (T^{1,0}_IN),\quad\mbox{for all}\quad 
X, Y \in \G (T^{1,0}_IN),\label{e:int}\\
d\tilde{\varphi}&=&0,\label{e:Kaehler}
\eea
where $\tilde{\varphi}$ is the pull back  of the fundamental form
$\varphi = h (\cdot ,J\cdot )$ to a two-form on $N$. Explicitly, 
\[ \tilde{\varphi}|_{\mathcal{D}\wedge TN}=0,\quad 
\tilde{\varphi}|_{\wedge^2\mathcal{N}}=\varphi.\]
The equation \re{e:isom} holds because the 
$\xi_i$ are Killing fields. The Lie derivative 
$\mathcal{L}_{\xi_i}I$ of $I\in \G (N,Q)$ is again a section
of $Q|_N$, since any isometry of a 
quaternionic K\"ahler manifold of nonzero scalar curvature
preserves the quaternionic structure $Q$. In order to prove 
\re{e:I}, it thus suffices to check that: 
\[ (\mathcal{L}_{\xi_i}I)\xi_1 = \mathcal{L}_{\xi_i}(I\xi_1)-
I\mathcal{L}_{\xi_i}\xi_1= \mathcal{L}_{\xi_i}(\xi_2)=0.\]
The equation \re{e:I} implies that on $N$ we have  
\[ P_i= (\n \xi_i)^Q 
= (\n_{\xi_i}-\mathcal{L}_{\xi_i})^Q \equiv \frac{1}{2}\o_1(\xi_i)J_1
+\frac{1}{2}\o_2(\xi_i)J_2\pmod{\bR J_3}.\]
Combining this with the equation $P_1P_2|_N=fJ_3|_N$ we obtain that, on $N$, 
\be P_i = \frac{1}{2}\o_1(\xi_i)J_1
+\frac{1}{2}\o_2(\xi_i)J_2,\ee
where the vectors $v_i := (\o_1(\xi_i),\o_2(\xi_i))\in \bR^2$ satisfy 
\be v_1\perp v_2,\quad |v_1||v_2|=4f
\ee on $N$. 
This shows that the one-forms $\o_1$, $\o_2$ 
are pointwise linearly independent on $N$. In fact, their restrictions 
to $\mathcal{D}$ are linearly independent. Hence,  
\[ \mathcal{K} := \ker \o_1|_N \cap  \ker \o_2|_N\subset TN\]
is a distribution complementary to $\mathcal{D}$. 
We will now show that  $\mathcal{K}=\mathcal{N}$. 
To see this, we calculate the covariant derivative of $P_i$ on $N$:
\bean \n P_i &=&  \frac{1}{2}(\n \o_1(\xi_i))\ot J_1
+\frac{1}{2}(\n \o_2(\xi_i))\ot J_2\\
&& -\frac{1}{2}\o_1(\xi_i)(\o_2\ot J_3 
-\o_3\ot J_2) 
-\frac{1}{2}\o_2(\xi_i)(\o_3 \ot J_1-\o_1\ot J_3)\\
&\equiv& -\frac{1}{2}\left( \o_1(\xi_i)\o_2 - \o_2(\xi_i)\o_1 \right) \ot J_3 \pmod{T^*M\ot (\bR J_1 \oplus \bR J_2)}.\eean
Comparing with \re{e:moment}, we obtain 
\[ \nu \varphi_3 (\xi_i,\cdot ) = \o_1(\xi_i)\o_2 - \o_2(\xi_i)\o_1\]
along $N$. This implies that $\mathcal{K}=\mathcal{N}$. It also shows that 
the two-form $\nu \varphi_3-\o_1\wedge\o_2$ vanishes on $\mathcal{D}\wedge 
TM$ along $N$ and coincides with $\nu \varphi_3$ on $\mathcal{N}$. This means that
\[ \nu \tilde{\varphi} = (\nu \varphi_3-\o_1\wedge\o_2 )|_N 
\stackrel{\re{e:str}}{=}
d\o_3|_N,\]
proving \re{e:Kaehler}. It remains to check the integrability
condition \re{e:int}, which shows that $I$ defines a 
complex structure on $N$. The equation \re{e:omega} implies that 
$\n_XI=0$, for all $X\in \mathcal{N}=\mathcal{K}$. Using the 
symmetry of the Levi-Civita connection and the fact that
$I\mathcal{N}=\mathcal{N}$, we can easily check that  
\bean [ X -iIX, Y - iIY] &=& [X,Y] - [IX,IY] -i([X,IY] + [IX,Y])\\
&=& [X,Y] - [IX,IY] -iI([X,Y] - [IX,IY]),\eean
for all $X, Y\in \G (\mathcal{N})$. It 
remains to calculate $[\xi_1-i\xi_2,Y - iIY]$ for any $Y\in \mathfrak{X}(N)$:
\[ [\xi_1-i\xi_2,Y - iIY] \stackrel{\re{e:I}}{=} 
[\xi_1-i\xi_2,Y] -iI[\xi_1-i\xi_2,Y]\in 
T^{1,0}_IM.\] 
This proves \re{e:int}.
\qed

\bc \label{cor:main} If, in addition to the assumptions of Theorem \ref{thm:main}, 
the vector fields $\xi_1|_N$, $\xi_2|_N$ generate a free and proper action of 
a two-dimensional Abelian Lie group $A$ on the submanifold $N\subset M$, 
then the quotient
$M':=N/A$ is a smooth manifold, which inherits a K\"ahler structure 
$(h, J)$ from the transversal geometry of the integrable distribution
$\mathcal{D}$. 
The projection $(N,g) \ra (M',h)$ is a Riemannian submersion and
a principal fiber bundle with structure group $A$. Moreover,
$(N,I) \ra (M',J)$ is holomorphic, where $I\in \G (N,Q)$ is 
the (integrable) almost complex structure which maps $\xi_1|_N$ to $\xi_2|_N$. 
If, more generally, the proper action of $A$ is only locally free with 
finite stabilisers, then $(M',h,J)$ is a K\"ahler orbifold. 
\ec
\section{Examples}\label{sect:ex}

\subsection{Hyperbolic 4-space}
As a first example, let us consider the four-dimensional hyperbolic space 
\begin{displaymath}
M = H_{\bR}^4=\frac{\SO_0(1,4)}{\SO(4)} \ .
\end{displaymath}
The solvable Iwasawa subgroup $L$ of $\SO_0(1,4)
=\mathrm{Isom}_0(M)$ acts simply 
transitively on $M$ and we can identify $M$ with the group manifold 
$L$ endowed with a left-invariant metric $g$ of constant curvature $-1$. 
$(M,g)$ is a quaternionic K\"ahler manifold with the  
quaternionic structure $Q$ spanned by three left-invariant 
complex structures $J_\a$, $\a=1,2,3$. 
The Lie algebra 
\[ \gl := \mathrm{Lie}\, L = \ga + \gn\] 
is the orthogonal sum of a 
three-dimensional Abelian nilradical $\gn = \mathrm{span} \{ X_\a = J_\a X_0
|\a=1,2,3\}$
and a one-dimensional subalgebra $\ga = \bR X_0$, where $X_0$ is a unit vector 
such that $ad_{X_0}|_\gn = \mathrm{Id}$.  Decomposing the Levi-Civita connection
\[ \n_XY = \frac{1}{2}\sum \o_\a (X)J_\a Y + \bar{\n}_X Y,\quad X,Y\in \gl ,\]
such that $\bar{\n}_XJ_\a=0$, one can easily compute $\o_\a=-X_\a^*$, where
$(X_a^*)$ is the dual basis of $\gl^*$.

Let us use  $k_a$, $a=0,1,2,3$, to denote
the (right-invariant) Killing vector field which coincides with 
the left-invariant vector field $X_a$ at $e\in L$. A straightforward calculation
shows that 
\begin{eqnarray*} k_0 (p) &=& X_0(p) - \frac{e^{-x^0}-1}{x^0}\sum x^\a X_\a(p)\\
k_{\a} &=& e^{-x^0}X_\a , 
\end{eqnarray*} 
at $p= \exp (x)\in L=M$, where $x=\sum x^\a X_\a$. 
This allows us to compute the moment maps $P_\a$ of the three commuting Killing 
vector fields $\k_\a$:
\[ P_1 = (\n k_1)^Q = -(\mathcal{L}_{k_1})^Q+\n_{k_1 }^Q = \frac{1}{2}\sum \o_\a (k_1)J_1 = 
-\frac{1}{2}e^{-x^0}J_1,\]
since $(\mathcal{L}_{k})^Q=0$ for any right-invariant Killing vector field $k$ 
and $\o_\a=-X_\a^*$. Summarising, we have shown that  
\begin{equation}
P_\alpha = -\frac{1}{2}e^{-x^0}J_\alpha \ , 
\end{equation}
in accordance with \cite{Ferrara:1995gu}. Thus we have   
\begin{equation}
P_1 P_2 k_1 = f k_2 \ ,\quad 4f=|k_1|^2=|k_2|^2=e^{-2x^0} >0, 
\end{equation}
and we can choose $\xi_i=k_i$, $i=1,2$, in agreement with conditions $(i-ii)$ 
in Theorem~\ref{thm:main}. The Killing vector fields $k_1, k_2$ generate
the left-action of the normal subgroup $A_2=\exp \ga_2 \subset L$, 
$\ga_2 =  \mathrm{span} \{ X_1, X_2\}$. 
Therefore, we can apply Theorem~\ref{thm:main} and 
Corollary~\ref{cor:main} to $N=M$. The quotient $M'=M/A_2$ is 
the complex hyperbolic line $M'\cong H_\bC^1$, which again has 
constant curvature $-1$ and admits the simply transitive solvable
group $L/A_2$ of holomorphic isometries. 

\subsection{Quaternionic K\"ahler manifolds in the image of the c-map}
There is a class of quaternionic K\"ahler manifolds of negative 
scalar curvature of the form $M=M_{\rm  sk}\times G$, where 
$M_{\rm  sk}$ is a (projective)  
special K\"ahler manifold 
of dimension $2n-2$ and $G$ is the solvable Iwasawa subgroup 
of $\SU(1,n+1)$, which is a semidirect product of a 
$(2n+1)$-dimensional Heisenberg group with $\mathbb{R}$.
For simplicity we will assume from now on that $M_{\rm  sk}$ 
admits a global system of special coordinates. Such manifolds
are called (projective)  {\cmssl special K\"ahler domains}. 
Note that the quaternionic K\"ahler metric on $M$ cannot be 
a product metric, 
since quaternionic K\"ahler manifolds are irreducible.  
The construction of these manifolds out of the special K\"ahler base 
is called the (supergravity) c-map \cite{Cecotti:1988qn,Ferrara:1989ik}. It has 
been recently shown that 
the quaternionic K\"ahler manifold $M$ is complete if $M_{\rm  sk}$ is 
complete \cite{CMX}.  As we will show, the class of quaternionic K\"ahler manifolds in the 
image of the c-map gives numerous examples for the quotient 
construction introduced in Theorem \ref{thm:main}.

In the following we will briefly describe the construction of the c-map, see
\cite{Cecotti:1988qn,Ferrara:1989ik,CMX} for more detailed information. 
Any  (projective) special  K\"ahler manifold $M_{\rm  sk}$ can be 
realised as the base  of a holomorphic $\bC^*$-principal bundle
$M_{\rm  ask}\ra M_{\rm  sk}$. The total space
$M_{\rm  ask}$ has the structure of an affine special
K\"ahler manifold, which admits special holomorphic local 
coordinates 
$Z^A$, $A=1,\ldots ,n$, such that the geometric data of $M_{\rm  ask}$
are encoded in a holomorphic function $F(Z^1,\ldots , Z^n)$ called
the {\cmssl holomorphic prepotential}.\footnote{Readers familiar with the 
supergravity literature might prefer to label the coordinates by 
$I=0,1,\ldots ,n-1$, as is done from Section \ref{sect:cubic} onwards.} 
The functions $z^a=Z^a/Z^n$, $a=1,\ldots ,n-1$, induce local coordinates 
on $M_{\rm  sk}$ and the K\"ahler potential $K(z)$ of $M_{\rm  sk}$ 
can be explicitly
expressed in terms of the prepotential $F$. In fact, $K(z)=K(z,1)$, where 
\be \label{e:Kpot} 
K(Z) = -\ln (2Z^AN_{AB}\bar{Z}^B),\quad N_{AB} = \Im F_{AB},\ee 
where the subscripts on $F$ denote derivatives with respect to $Z^A$ e.g. ${F}_A=\partial{F}/\partial Z^A$. 
The solvable Lie group $G$ admits
a natural system of global coordinates
\footnote{In supergravity theories arising as effective theories of type II
compactifications the scalar manifold $M_{\rm  sk}$ is spanned by deformations of the metric and the Neveu-Schwarz B-field, while $(\phi , \tilde \phi , a^A, b_A)$  correspond to the dilaton, the axion and the 2n real Ramond-Ramond scalars, respectively.}
$(\phi , \tilde \phi , a^A, b_A)$, $A=1,\ldots,n$. 
A basis for the right-invariant vector fields on $G$ 
is given in these coordinates by 
\begin{equation}\label{e:Killing}
\begin{aligned}
k_{\phi}   ~= &~ \tfrac{1}{2} \frac{\partial}{\partial \phi} -  {\tilde \phi} \frac{\partial}{\partial {\tilde \phi}} - \tfrac{1}{2} a^A \frac{\partial}{\partial a^A} - \tfrac{1}{2} b_A \frac{\partial}{\partial b_A} \ , \\
k_{\tilde \phi}   ~= &~ - 2 \frac{\partial}{\partial {\tilde \phi}} \ , \\
k_A  ~= &~ \frac{\partial}{\partial a^A} + b_A \frac{\partial}{\partial {\tilde \phi}} \ , \\
\tilde k^A  ~= &~ \frac{\partial}{\partial b_A} - a^A \frac{\partial}{\partial {\tilde \phi}} \ .
\end{aligned}
\end{equation}
These vector fields obey the commutation relations
\begin{equation}\label{e:isometries_fibre}
\begin{aligned}
&& [k_{\phi},k_{\tilde \phi}]\  =& k_{\tilde \phi} \ , \qquad \qquad
&& [k_{\phi},k_A]\ \ = & \tfrac{1}{2} k_A \ ,  \\
&& [k_{\phi},\tilde k^A]\  =& \tfrac{1}{2} \tilde k^A \ , \qquad \qquad
&& [k_A,\tilde k^B]\  = &  \delta_A^B k_{\tilde \phi} \ ,
\end{aligned}
\end{equation}
with all other commutators vanishing.

Recall that a {\cmssl quaternionic vielbein} on a 
quaternionic K\"ahler manifold  (or, more generally, on an almost
quaternionic Hermitian manifold) $(M,g,Q)$ is a coframe which belongs to the 
$\Sp (n)\Sp(1)$-structure defined by $(g,Q)$, cf.\ \cite{ACDGV}. 
More explicitly, it is a system of complex-valued one-forms 
$\mathcal U^{\mathcal A m}$, 
$\mathcal A = 1,2, m=1,\dots,2n,$ such that the metric takes the form
\be \label{gEqu} g = \sum \e_{\mathcal{A}\mathcal{B}}\e_{lm}
\mathcal U^{\mathcal A l}\ot \mathcal U^{\mathcal B m}\ee
and such that the quaternionic structure $Q$ on $TM$ corresponds
to the standard quaternionic structure on the first factor $\bC^2$
of the tensor product $\bC^2\ot \bC^{2n}$. Here 
$\e = \left( 
\begin{array}{cc} 0& \mathbbm{1}\\
-\mathbbm{1} & 0
\end{array}\right)$. Note that the 
metric and quaternionic structure are completely determined
by specifying a quaternionic vielbein. 

In \cite{Ferrara:1989ik} it was proven that 
\begin{equation} \label{e:quat_vielbein}
\mathcal U^{\mathcal A m}= \tfrac{1}{\sqrt{2}}
\left(\begin{aligned}
\bar{u} && \bar{e} && -v && -E \\
\bar{v} && \bar{E} && u && e
\end{aligned}\right)
\end{equation}
is a quaternionic vielbein of a  
quaternionic K\"ahler structure 
$(g,Q)$ on a domain $M$ 
if the one-forms $\mathcal U^{\mathcal A m}$ are defined by
\begin{equation} \label{e:one-forms_quat}
\begin{aligned}
u ~= &~ \iu e^{K/2+\phi}Z^A(\diff b_A - {F}_{AB} \diff a^B) \ , \\
v ~= &~ \tfrac{1}{2} e^{2\phi}\big[ \diff e^{-2\phi}-\iu (\diff {\tilde \phi} 
+b_A \diff a^A-a^A \diff b_A  ) \big] \ , \\
E^{\,b} ~= &~ -\tfrac{\iu}{2} e^{\phi-K/2} {\Pi}_A^{\phantom{A}b} 
N^{AB}(\diff b_B - {\bar {F}}_{BC} \diff a^C) \ , \\
e^{\,b} ~= &~ {\Pi}_A^{\phantom{A}b} \diff Z^A = 
e^{\,b}_adz^a\ . \\
\end{aligned}
\end{equation}
Here $(Z^A)$, $A=1,\ldots, n$, 
are the homogeneous coordinates of $M_{\rm  sk}$, 
which are functions on the affine special K\"ahler domain $M_{\rm  ask}$,  
\[{\Pi}_A^{\phantom{A}b} =(e_a^{\phantom{a}b},
-z^ae_a^{\phantom{a}b})\] 
is defined 
using the vielbein $e_a^{\phantom{a}b}$ on 
$M_{\rm  sk}$. 
In the above formulas
one may simply put $(Z^A)=(z^a,1)$ to obtain differential forms
which are manifestly defined on $M=M_{\rm sk}\times G$, rather than horizontal 
$\bC^*$-invariant forms on $M_{\rm ask}\times G\ra M_{\rm sk}\times G$. 
It is shown in \cite{CMX} that, although the prepotential $F$ and the 
vielbeins are coordinate dependent, 
the quaternionic K\"ahler structure does not depend (up to isomorphism) 
on the choice of special coordinates. 

\noindent 
{\bf Remark:} 
It is also shown in \cite{CMX} that a global quaternionic K\"ahler structure 
can be defined even if $M_{\rm sk}$ cannot be covered by a single 
system of special coordinates. In that case one has to replace 
$M=M_{\rm sk}\times G$  
by the total space of a possibly nontrivial bundle over $M_{\rm sk}$. 

Using the quaternionic vielbein given in \eqref{e:quat_vielbein} we can define 
three almost complex structures $J_\alpha$ on $M$ by 
\begin{equation} \label{e:conn_J}
\mathcal U^{\mathcal A m}\circ J_\alpha  = - \iu (\sigma_\alpha)^{\cal A}_{~\cal B} \mathcal U^{\mathcal B m} \ ,
\end{equation}
where $(\sigma_\alpha)^{\cal A}_{~\cal B}$ are the $su(2)$ generators
\begin{displaymath}
(\sigma_1)^{\cal A}_{~\cal B} = \left( \begin{aligned} 0 && 1 \\ 1 && 0\end{aligned} \right) \ , \qquad (\sigma_2)^{\cal A}_{~\cal B} = \left( \begin{aligned} 0 && -\iu \\ \iu && 0\end{aligned} \right) \ , \qquad (\sigma_3)^{\cal A}_{~\cal B} = \left( \begin{aligned} 1 && 0 \\ 0 && -1\end{aligned} \right) \ . 
\end{displaymath}
Then $Q =\mathrm{span}\{ J_\a| \a = 1,2,3\}$ is a   
skew-symmetric parallel quaternionic structure with respect to the 
quaternionic K\"ahler metric 
\be \label{e:metric_cmap}
 g=u\bar u+ v\bar v +  
\sum (e^{\,b} \bar e^{\,b} + 
E^{\,b} \bar E^{\,b}) \ee 
on $M$ defined by
\re{gEqu}. (Recall that $u\bar u =\bar u u= \frac{1}{2}(u\ot \bar u + 
\bar u \ot u)$.)  

\bp \label{J3Prop} The almost complex structure $J_3$ is integrable for 
any quaternionic K\"ahler manifold $(M=M_{\rm sk}\times G,g,Q)$ in the 
image of the c-map. Moreover, the factors of the product $M_{\rm sk}\times G$
are complex submanifolds of the complex manifold $(M,J_3)$. The 
restriction of $J_3$ to the first factor coincides (at any point of $M$) with 
the original complex structure $J$ on the K\"ahler manifold $M_{\rm sk}$, 
whereas the submanifold $G= \{ p\}\times G\subset M_{\rm sk}\times G =M$ 
with the Hermitian structure induced by $(g,J_3)$ is not K\"ahler.  
Nevertheless, the submanifold $G\subset M$ with its induced metric 
is isometric to the complex hyperbolic space 
$H_\bC^{n+1}$ with the K\"ahler metric of constant holomorphic sectional 
curvature $-4$. \ep

\pf
According to \eqref{e:quat_vielbein} and \eqref{e:conn_J}, the 
one-forms $u, \bar v, e^{\,b}, \bar E^{\,b}$ 
constitute a  basis for the
space of $(1,0)$-forms of $J_3$. We can compute their exterior 
derivative to be \cite{Ferrara:1989ik} 
\begin{displaymath}
\begin{aligned}
\diff u =& \left( -\tfrac12 (v+ \bar v) + \frac{\bar Z^A N_{AB}\diff Z^B - Z^A N_{AB}\diff \bar Z^B}{2\bar Z^A N_{AB} Z^B} \right)\wedge u - \bar E \wedge e \ , \\ 
\diff \bar v =& - v \wedge \bar v +\bar u \wedge u - E \wedge \bar E \ , \\
\diff e^{a} =& - \omega^{a}_{~b} \wedge e^{b} \ , \\
\diff \bar E^{a} =& \left( -\bar \omega^{a}_{~b} - \tfrac12 (v+\bar v)\d^{a}_{b} + \frac{\bar Z^A N_{AB}\diff Z^B - Z^A N_{AB}\diff \bar Z^B}{2\bar Z^A N_{AB} Z^B}\d^{a}_{b} \right) \wedge \bar E^{b}\\ & + \bar e \wedge u + \tfrac14 (\bar Z^A N_{AB} Z^B)   \Pi_{b\,A} N^{AB} N^{CD} \Pi^{a}_D E^{b} \wedge \diff {F}_{BC}   \ ,
\end{aligned}
\end{displaymath}
where $\omega$ is the connection one-form of $M_{\rm  sk}$ and the index 
$b$ is lowered by means of the Kronecker symbol. Since there is no $(0,2)$-form appearing on the right-hand side, $J_3$ is integrable in virtue
of the Newlander-Nirenberg theorem. The two distributions tangent to the 
factors of the product manifold $M_{\rm sk}\times G$  are defined by
$u=v=E^{\,b}=\bar u=\bar v=\bar E^{\,b}=0$ and 
$e^{\,b}=\bar e^{\,b}=0$, respectively. This shows
that both distributions are $J_3$-invariant and, hence, that
the leaves are complex submanifolds. The formula $e^{\,b}=
e^{\,b}_adz^a$ implies that the complex structures $J_3|_{M_{\rm  sk}}$
and $J$ coincide. It is known that a K\"ahler submanifold $S\subset M$ of a 
quaternionic K\"ahler manifold $M$, such that the complex structure of $S$ 
is subordinate to the quaternionic structure,  
has at most dimension $\frac{1}{2}\dim M$ \cite{AM}. Since $\dim G = 2n+2$, 
$2n=\frac{1}{2}\dim M$, $G\subset M$ cannot be a K\"ahler submanifold with
the complex structure induced by $J_3$. Alternatively, one may check by a 
direct calculation that 
the fundamental two-form $\varphi_3 = g(\cdot , J_3)$ is not closed. 
For a proof of the last statement of the 
proposition see \cite{CMX}. 
\qed

In the next proposition we give more detailed information about
the complex structure $J_3$. 
\bp  \label{biholProp} \begin{enumerate}
\item[(i)] The complex structure $J_3$ on the quaternionic K\"ahler manifold 
$M=M_{\rm sk}\times G$
is of the form $J_3 = J + J_G$, where $J$ is the complex structure of
the projective special K\"ahler domain $M_{\rm sk}$ and 
$(J_G(p))_{p\in M_{\rm sk}}$ is a smooth family of left-invariant
complex structures $J_G(p)$ on $G$. 
\item[(ii)] The projection 
$\pi : M \ra M_{\rm sk}$ is a holomorphic submersion with  
fibers $(G,J_G(p))$ biholomorphic to the domain
\[ F(n+1) := \{ (w^0,w_1,\ldots ,w_n)\in \bC^{n+1}| \Re w^0 > 
\sum_{A=1}^{n-1} (\Im w_A)^2-(\Im w_n)^2
\} \subset \bC^{n+1},\]
for all $p\in M_{\rm sk}$.  
The total space $(M,J_3)$ admits a fiber preserving open  
holomorphic embedding 
into the trivial holomorphic bundle $M_{\rm sk}\times \bC^{n+1}$. 
\end{enumerate}
\ep 
\pf (i) By Proposition \ref{J3Prop}, the complex structure $J_3$ on 
$M=M_{\rm sk}\times G$ is the sum 
of the complex structure $J$ on the base and a family of complex structures
$J_G(p)$ on the fibers $\{ p\} \times G \cong G$. To prove that 
$J_G(p)$ is left-invariant it suffices to check that the Lie derivative
of the one-forms \re{e:one-forms_quat} with respect to the right-invariant 
vector fields  \re{e:Killing} vanishes. That is a straightforward 
calculation.\\  
(ii) We define a fiber preserving holomorphic embedding $\Psi : 
M \ra M_{\rm sk}\times
\bC^{n+1}$ by $\Psi = (\pi , w^0,w_A)$, where 
\be \label{holcoordEqu} 
w^0:=e^{-2\phi} + i(\tilde \phi +a^A(b_A -F_{AB}a^B)),\quad 
w_A:= b_A-F_{AB}a^B.\ee 
One can easily check that the functions $w^0, w_A$ are 
$J_3$-holomorphic, cf.\ \cite{LST2}.  We claim that $\Psi$ maps
$M$ biholomorphically onto the domain defined by the 
inequality 
\be \label{ineq} \Re w^0 > -N^{AB}
\Im w_A \Im w_B. 
\ee 
In fact, for fixed $p\in M_{\rm sk}$, the linear map 
\[ \bR^{2n}\ni (a^A,b_B) 
\mapsto (w_A)\in \bC^n\] 
is an isomorphism, whereas the variable  
$w^0=e^{-2\phi} + i(\tilde \phi +a^Aw_A)$ is constrained by the 
inequality 
$\Re w^0 >  a^A\Im w_A$.
Expressing $(a^A)$ by $(w_A)$
yields 
\[ a^A = -N^{AB}\Im w_B\]
and thus \re{ineq}. 
For fixed $p\in M_{\rm sk}$ we can choose the special coordinates
such that $(N^{AB}(p))=\mathrm{diag}(-1,\ldots, -1,1)$. This shows that 
$\pi$ is a holomorphic submersion with fibers biholomorphic to $F(n+1)$. 
\qed  

Given the explicit form of the vielbein \eqref{e:quat_vielbein} the $\SU(2)$ connection $\omega^x$ reads \cite{Ferrara:1989ik}
\begin{equation}\label{e:quat_connection}
\begin{aligned}
\omega^1 ~= &~  \iu (\bar u- u)  \ , \qquad
\omega^2 ~= ~ u + \bar u \ , \\
\omega^3 ~= &~ \tfrac{\iu}{2} (v-\bar v) - \iu e^{K} \left(Z^A N_{AB} \diff \bar Z^B - \bar Z^A N_{AB} \diff Z^B \right) \ .
\end{aligned}
\end{equation}
It can be checked that the natural action of $G$ on $M=M_{\rm sk}\times G$
preserves the Ferrara-Sabharwal metric $g$ \cite{CMX}. 
The moment maps $P_\lambda$ of the 
Killing vectors $k_\lambda$ given in \eqref{e:Killing} take the 
following simple form 
\begin{equation} \label{e:prepotential_no_compensator}
P_\lambda = \frac{1}{2}\sum \omega_{\alpha} (k_\lambda) J_\alpha \ .
\end{equation}
This follows from $\n k_\l = \n_{k_\l}-\mathcal{L}_{k_\l}$, since 
$\mathcal{L}_{k_\l}J_\a=0$.\footnote{Note that the formula \re{e:prepotential_no_compensator} 
differs by a factor $1/2$ from that of \cite{Michelson:1996pn},
since our definition \re{momentEqu} of the moment map differs from that
of \cite{Michelson:1996pn} by the same factor. This can be easily checked
with the help of formula \re{e:moment}.} 

In order to define the submanifold $N\subset M$ to which we will apply 
the quotient construction of Theorem~\ref{thm:main}, we choose constant complex vectors $(C_A)$ and $(D^A)\neq 0$ and a constant $\tilde C$, 
where $D^A$ 
obeys
\begin{equation}\label{e:null}
\sum_{A, B=1}^n  N_{AB}(Z_0) D^A \bar D^B = 0 \ 
\end{equation} 
at some point $Z_0 = (Z^1_0,\ldots ,Z^n_0)\in M_{\rm ask}$. Here 
$M_{\rm ask}$ is identified with a domain in $\bC^n$ by means of the special
coordinates. Since the affine special K\"ahler 
metric $\sum N_{AB}(Z_0)dZ^Ad\bar{Z}^B$ 
is indefinite, such a vector $(D^A)$ does always exist. 
We will assume that 
the rank of the matrix 
\be \label{equ:GAB} G_{AB}(Z) = \sum_C D^CF_{ABC}(Z)\ee 
is constant 
in a neighborhood of $Z_0$. Then, by restricting to that neighborhood, we
can assume that the rank is constant on $M_{\rm ask}$. 
That implies that the system of equations 
\be \label{MwedgeEqu} \sum_A D^A {F}_{AB}(Z) = C_B := 
\sum_A D^A {F}_{AB}(Z_0)\ee 
defines a complex submanifold $M^\wedge_{\rm ask}\subset M_{\rm ask}$ of 
complex dimension $n-r$, where $r = \mathrm{rk} (G_{AB})$. 
\bp \label{rankProp} $r\le n-1$ and $M^\wedge_{\rm ask}$ fibers over a 
complex submanifold $M^\wedge_{\rm sk}\subset M_{\rm sk}$ of dimension
$n-1-r$. In particular,  
$M_{\rm sk}$ is of dimension zero if $(G_{AB})$ has maximal rank.
\ep 
\pf Since $F_{AB}$ is homogeneous of degree zero, the vector
$(Z^A)$ is in the kernel of the matrix $(G_ {AB})$, which implies $r\le n-1$. 
Due to the 
homogeneity of the equation \re{MwedgeEqu}, 
$M^\wedge_{\rm ask}$ is a cone over a complex submanifold 
$M^\wedge_{\rm sk}\subset M_{\rm sk}$. \qed 
 
\noindent
{\bf Remark:} \label{GAB}
More generally, for the smoothness of $M^\wedge_{\rm ask}$ it is 
sufficient to assume
that the rank of $(G_{AB})$ is constant on a complex submanifold containing
(a neighborhood of $Z_0$ in) the analytic set defined by \re{MwedgeEqu}. 

Now we define a subset $N\subset M$ by the system 
\begin{equation} \label{e:definition_N}
D^A {F}_{AB}(Z) = C_B \ , \qquad
D^A (b_A - {F}_{AB} a^B ) = \tilde C \ .
\end{equation}
We claim that $N\subset M$ is submanifold of codimension $2r+2$.  
More precisely, it is a subbundle of $M^\wedge_{\rm sk}\times G$
with fibers of codimension 2. To see this it suffices to recall 
that the first equation  of \re{e:definition_N} defines the 
submanifold $M^\wedge_{\rm sk}\subset M_{\rm sk}$ and to note that over points of 
$M^\wedge_{\rm sk}$ the second equation reduces to 
\[ D^A b_A - C_B a^B  = \tilde C,\]
which is a system of two real affine
equations. The two real equations are independent
if and only if the vector $(D^A,C_B)$ is not a complex multiple of a 
real vector. The latter property follows from the fact that that vector
belongs to the tangent space $L=T_{dF(Z_0)}\mathcal{C}$ of the Lagrangian cone 
$\mathcal{C}=\{dF(Z)|Z\in M_{\rm ask}\subset \bC^{n}\}\subset T^*\bC^n\cong \bC^{2n}$, 
which satisfies $L\cap \bar{L}=0$, see \cite{ACD}.  

\bp \label{NProp} Under the above assumptions, 
$N\subset M$ is a complex submanifold with respect to
the complex structure $J_3$. More precisely, the homogeneous equation
\be \label{homEqu} D^A b_A - C_A a^A  = 0\ee
defines a subgroup $G^\wedge \subset G$ of codimension 2 and 
$N =M^\wedge_{\rm sk}\times S$ is the product of the 
complex submanifold $M^\wedge_{\rm sk}\subset M_{\rm sk}$ and a  
submanifold $S\subset G$, which is a left-translate,   
$S=xG^\wedge$,  of the subgroup $G^\wedge\subset G$ by an element  
$x\in G$ satisfying the inhomogeneous equation 
\be  \label{inhomEqu} D^A b_A - C_A a^A  = \tilde C.\ee  
The fibers $\{ p\} \times S \subset \{ p\} \times G$, 
$p\in  M^\wedge_{\rm sk}$, 
are complex hypersurfaces with respect to the 
complex structure on $\{ p\} \times G\subset M$ induced by $J_3$. 
\ep 

\pf In order to prove that $N\subset M$ and the fibers of $N\ra  
M^\wedge_{\rm sk}$ are a complex submanifolds,  
it suffices to show that the one-form 
\[ d(D^A(b_A-F_{AB}a^B))= D^A(db_A -F_{AB}da^B) + D^Aa^BdF_{AB}\]
is of type $(1,0)$. This is obvious for the second term.  
In order to analyse the first term, we decompose
\[ D^A = cZ^A + H^A,\]
where $\bar{H}^AN_{AB}Z^B=0$ and $c\in \bC$. Then 
\[ D^A(db_A -F_{AB}da^B) = -cie^{-K/2-\phi} u +
c_{\,b}\bar{E}^{\,b},\]
where the coefficients $c_{\,b}\in \bC$ are determined by the equation 
$-\tfrac{\iu \bar{c}_{\,b} }{2} 
e^{\phi-K/2}{\Pi}_A^{\phantom{A}b} N^{AB}=\bar{H}^B$. This proves
that $D^A(db_A -F_{AB}da^B)$ is of type $(1,0)$. 

To check that \re{homEqu} defines a subgroup $G^\wedge \subset G$, we 
recall\footnote{Our additive variable $\phi$ is related to the   
corresponding variable $\l$ in \cite{CMX} by $\l = -2\phi$.} 
\cite{CMX} that in the coordinates $(\phi , \tilde \phi , a^A, b_B)$ the group 
multiplication in $G$ is given by:
\be \label{multEqu} (\phi , \tilde \phi , a, b)\cdot (\phi' , \tilde \phi' , a', b')
= (\phi + \phi',  \tilde \phi + e^{-2\phi} \tilde \phi' + 
e^{-\phi}(a^Ab'_A-a'^Ab_A), a+e^{-\phi}a',b+e^{-\phi}b').
\ee 
{}From this formula we see that the set of solutions of \re{homEqu} is closed
under multiplication and contains the neutral element and the 
inverse 
\[ (\phi , \tilde \phi , a, b)^{-1}=(-\phi , -e^{2\phi}\tilde \phi,
-e^{\phi}a,-e^{\phi}b)\]
of any element $(\phi , \tilde \phi , a, b)$ satisfying \re{homEqu}.   
Let $x\in G$ be any element satisfying \re{inhomEqu}. Using the 
multiplication law \re{multEqu} we can easily check that 
$xG^\wedge$ is a subset of the solution space of \re{inhomEqu}, which we
know is an affine subspace of $G\cong \bR^{2n+2}$ of codimension 2. This 
proves that $xG^\wedge$ coincides with the set of solutions of \re{inhomEqu},
that is, with the fiber of $N\ra M^\wedge_{\rm sk}$. 
\qed 

In the next proposition we give more detailed information about
the complex submanifold $N\subset (M,J_3)$. 
\bp  \label{NbiholProp} \begin{enumerate}
\item[(i)] The complex structure induced by $J_3$ 
on $N=M_{\rm sk}^\wedge\times S$ is of the form 
$J+J_S$, where $J$ is the complex structure
on $M_{\rm sk}^\wedge$ and $(J_S(p))_{p\in M_{\rm sk}^\wedge}$ is a smooth family of 
left-invariant complex structures 
on $S=xG^\wedge\cong G^\wedge$. 
\item[(ii)] The projection 
$\pi_N: N \ra M_{\rm sk}^\wedge$ is a holomorphic submersion with 
fibers  $(S,J_S(p))$ biholomorphic to $B^{n-1}_{\bC}\times \bC$, 
for all $p\in M_{\rm sk}^\wedge$.  
The total space 
$(N,J_3)$ admits a fiber preserving open holomorphic embedding
into the trivial holomorphic bundle 
$M_{\rm sk}^\wedge \times \bC^n$. 
\end{enumerate}
\ep 

\pf (i) It follows from Proposition \ref{biholProp} 
and Proposition \ref{NProp} that the complex structure of $N$
is of the form $J+J_S$, where $J_S=J_G|_{S}$. Identifying 
$S=xG^\wedge$ with the group $G^\wedge\subset G$ by means of the 
left-translation with 
$x^{-1}$, we can consider $J_S$ as a complex structure on the 
group $G^\wedge$. Then the left-invariance
of $J_S$ follows from that of $J_G$.\\
(ii) Using the fiber preserving open holomorphic embedding $\Psi$ of 
$\pi : M \ra M_{\rm sk}$ into $M_{\rm sk} \times \bC^{n+1}$ defined in 
\re{holcoordEqu}, we see that $\pi_N: N \ra M_{\rm sk}^\wedge$  
is embedded into $\pi|_{M_{\rm sk}^\wedge}$ by one complex affine equation 
$D^Aw_A=\tilde{C}$, which reduces the 
trivial bundle $M_{\rm sk}^\wedge 
\times \bC^{n+1} \subset M_{\rm sk}\times \bC^{n+1}$ 
to a trivial bundle $\cong M_{\rm sk}^\wedge \times \bC^n$   
and the fiber $F(n+1)$ of $\pi$ to 
$F'(n-1)\times \bC$. In fact, for fixed $p\in M_{\rm sk}^\wedge$ we can choose
special coordinates such that  $(N^{AB}(p))=\mathrm{diag}(-1,\ldots, -1,1)$
and $D=(0,\ldots,0,1,1)$. Then the fiber is defined by 
\[ \Re w^0 > 
\sum_{A=1}^{n-1} (\Im w_A)^2 - (\Im w_n)^2,\quad w_{n}+w_{n-1}=\tilde{C}.\]
Elimination of $w_n$ yields the domain 
\[ \{ (w^0,\ldots ,w_{n-1})\in \bC^{n}| \Re w^0 - 2\Im \tilde{C}\Im w_{n-1}
+(\Im \tilde{C})^2> 
\sum_{A=1}^{n-2} (\Im w_A)^2\} \subset \bC^{n},\]
which is biholomorphic to $F'(n-1)\times \bC$ by the affine
transformation $(w^0,w_1,\ldots , w_{n-1}) \mapsto 
(w^0 - 2\Im \tilde{C}\Im w_{n-1}
+(\Im \tilde{C})^2,w_1,\ldots , w_{n-1})$, 
where 
\[ F'(n-1) := \{ (w^0,w_1,\ldots ,w_{n-2})\in \bC^{n}| \Re w^0 > 
\sum_{A=1}^{n-2} (\Im w_A)^2\} \subset \bC^{n-1}.\]
Now it suffices to note that $F'(n-1)$ is biholomorphic to
the ball $B_\bC^{n-1}$. 
\qed 

Before we go on, let us summarize what we found so far by the 
following commutative diagram consisting of holomorphic 
fiber preserving embeddings: 
\[ \begin{array}{ccc} M=M_{\rm sk}\times G &\xhookrightarrow{\Psi} &M_{\rm sk}\times \bC^{n+1}\\
\cup &&\cup\\
N= M_{\rm sk}^\wedge \times S  &\xhookrightarrow{\Psi|_N} &M_{\rm sk}^\wedge\times 
\bC^n, 
\end{array}
\]
where the horizontal embeddings are open and the vertical
ones are of complex codimension $r+1$. Recall that 
$r$ is the complex codimension of $M_{\rm sk}^\wedge \subset M_{\rm sk}$, 
$S$ is a left-translate of a subgroup $G^\wedge\subset G$ and
$\bC^n\subset \bC^{n+1}$ is an affine hyperplane (which is linear
if $S=G^\wedge$). The fibers of $M\ra M_{\rm sk}$ are biholomorphic
to $F(n+1)$, whereas the fibers of $N\ra M_{\rm sk}^\wedge$ are 
biholomorphic to $B^{n-1}_{\bC}\times \bC$.\\
 
Let us now define two Killing vectors $\xi_i$, $i=1,2$, on $M$ by
\begin{equation} \label{e:Killingv} \begin{aligned}
  \xi_1 =& \ \Re D^A k_A + \Re C_A \tilde k^A +\Re \tilde C k_{\tilde \phi}  \ , \\
  \xi_2 =&\ \Im D^A k_A + \Im C_A \tilde k^A +\Im \tilde C k_{\tilde \phi} \ .
\end{aligned}\end{equation}
From \eqref{e:prepotential_no_compensator}, \eqref{e:quat_connection}, \eqref{e:one-forms_quat} and \eqref{e:definition_N} we see that both $P_1$ and $P_2$ lie in the plane spanned by $J_1$ and $J_2$. Therefore, we find $P_1 P_2=f J_3$ for some function $f$. Furthermore, there is $J_3 k_1 = k_2$. 
Hence we can apply Theorem~\ref{thm:main} and Corollary~\ref{cor:main}, 
provided that $\xi_1, \xi_2$ are tangent to $N$ and generate a 
free and proper action. This is shown in the next proposition. 
\bp\label{holactionProp} 
The vector fields $\xi_1, \xi_2$ generate a free and proper
holomorphic action of a vector group $\bR^2\cong \bC$ on the submanifold 
$N\subset M$.   In the coordinates 
$(z^a,\phi,\tilde{\phi},a^A,b_A)$ the action of $(\l_1 ,\l_2)\in \bR^2$ 
is given by  
$(z,\phi,\tilde{\phi},a,b)\mapsto 
(z,\phi,\tilde{\phi}',a',b')$, 
where
\begin{eqnarray}\label{e:real} 
\tilde{\phi}' &=& \tilde{\phi}-\l_1\Re \tilde{C} -\l_2\Im \tilde{C},\\\nonumber
{a^A}' &=& a^A+\l_1\Re D^A + \l_2\Im D^A,\\\nonumber
b_A'  &=&b_A+\l_1\Re C_A + \l_2\Im C_A .
\end{eqnarray} 
In the holomorphic coordinates $(z^a,w^0,w_A)$
the action of $\l= \l_1 +i\l_2\in \bC$ 
is given by $(z^a,w^0,w_A)\mapsto (z^a,\zeta^0,\zeta_A)$,
where 
\begin{eqnarray}\label{e:holo}
\zeta^0 &=& w^0 +i\l \bar{D}^Aw_A -i\l \overline{\tilde{C}} + 
i\frac{\l^2}{4} \bar{D}^A(\bar{C}_A- F_{AB}\bar{D}^B),\\ \label{e:holo2}
\zeta_A &=& w_A + \frac{\l}{2} (\bar{C}_A-F_{AB}\bar{D}^B).
\end{eqnarray}
\ep  

\pf First Note that  
\begin{equation} \label{e:Killingvnew} 
\begin{aligned}
  \xi_1|_N =& \ \Re D^A \frac{\partial}{\partial a^A} + \Re C_A 
\frac{\partial}{\partial b_A} -\Re \tilde C \frac{\partial}{\partial \phi}  \ , \\
  \xi_2|_N =&\ \Im D^A \frac{\partial}{\partial a^A} + \Im C_A \frac{\partial}{\partial b_A} -\Im \tilde C \frac{\partial}{\partial \phi} \ .
\end{aligned}\end{equation}
We can easily check that $\xi_1 ,\xi_2$ are tangent to 
$N=\{ D^Ab_A-C_Aa^A=\tilde{C}\}$. In fact, this is a consequence of the
two equations $D^AC_A-C_AD^A=0$ and $D^A\bar{C}_A-C_A\bar{D}^A=-2i D^AN_{AB}
\bar{D}^B=0$.  Let us use $\varphi_j^t$ to denote the flow of the vector field 
$\xi_j$ and put 
\[ \varphi^\l := \varphi_1^{\l_1}\circ \varphi_2^{\l_2},\quad \l = \l_1+i\l_2.\]
Then \re{e:Killingvnew} shows that $\varphi^\l|_N$ is given by 
\re{e:real}. We see that in these coordinates the action consists of 
translations along a plane.  In particular,  
it is free and proper.  Expressing $\varphi^\l$ in holomorphic
coordinates yields \re{e:holo}-\re{e:holo2}, which shows that the action is 
$\bC \times N \ra N$ is holomorphic.
\qed 

The K\"ahler manifold $M'=N/A$ constructed from 
Corollary~\ref{cor:main} is of real dimension $4(n-1)-2r$, where 
$r$ was the complex codimension of $M^\wedge_{\rm sk}\subset M_{\rm sk}$.
Thus the minimal dimension of $M'$ is $2(n-1)$, which is attained when the 
base  manifold $M^\wedge_{\rm sk}$ is discrete. The maximal
dimension $4(n-1)$ is attained, when  $M^\wedge_{\rm sk}=M_{\rm sk}$. 
\bt \label{ThmM'} Let $(M',h)$ be the K\"ahler manifold 
obtained as above from the 
quotient construction  of Corollary~\ref{cor:main} applied to 
a quaternionic K\"ahler manifold $(M= M_{\rm sk} \times G,g)$ in the 
image of the c-map. Then $M'$ is the total space of a holomorphic 
submersion over the complex submanifold $M_{\rm sk}^\wedge \subset M_{\rm sk}$ with 
fibers biholomorphic to $B^{n-1}_{\bC}$.
The metric of the fiber is given by
\begin{equation}
\label{e:metric_fibre}
\begin{aligned}
h_{\rm fib} = & \tfrac{1}{4} e^{4 \phi} |\diff x^0 +2 \iu ( (\Im x)_a \d^{ab}) \diff x_b|^2 + \tfrac{1}{2} e^{2\phi} \diff \bar x_a \delta^{ab} \diff x_b \ ,
\end{aligned}
\end{equation}
with respect to some global system of holomorphic coordinates 
$(x^0,x_1,\ldots , x_{n-2})$ on the fiber. As a consequence, the 
fiber is isometric (but not biholomorphic, unless $n\le 2$) to $H_\bC^{n-1}$ 
with its metric of constant holomorphic sectional curvature $-4$. 
\et 

\pf Since the action on $N$ 
generated by $\xi_1$ and $\x_2$ 
is holomorphic, see Proposition \ref{holactionProp}, 
and preserves the fibers of 
the holomorphic submersion $\pi_N : N \ra M_{\rm sk}^\wedge$, 
we have an induced holomorphic submersion $\pi' : M'\ra M_{\rm sk}^\wedge$. 
We know already (see the proof of Proposition \ref{NbiholProp}) 
that $\pi_N: N\ra M_{\rm sk}^\wedge$ is holomorphically 
embedded into $\pi|_{M_{\rm sk}^\wedge} : \pi^{-1}(M_{\rm sk}^\wedge) \ra M_{\rm sk}^\wedge$ 
with fibers of complex codimension one and 
that $\pi^{-1}(M_{\rm sk}^\wedge)$ is an open subset of the 
trivial bundle $M_{\rm sk}^\wedge\times \bC^{n+1}$. The fiber $S=xG^\wedge$ of 
$\pi_N$ is the intersection of the fiber $G\cong F(n+1)$ of 
$\pi|_{M_{\rm sk}^\wedge}$ with the complex affine hyperplane 
defined by the equation $D^Aw_A=\tilde{C}$ in the holomorphic fiber
coordinates $(w^0,w_A)$. Let $V^A$ be any vector such that
$V^A(\bar{C}_A-F_{AB}\bar{D}^B)\neq 0$ holds on some neighborhood
$U\subset M_{\rm sk}^\wedge$. Such a vector exists, since  
$\bar{C}_A-F_{AB}\bar{D}^B= (\bar{F}_{AB}-F_{AB})\bar{D}^B=-2iN_{AB}\bar{D}^B$
and $D\neq 0$. Consider the subgroup $G'\subset G$ defined
by the homogeneous equations $D^Aw_A=V^Aw_A=0$. One can check that 
$G'$ is isomorphic to the Iwasawa subgroup
of $\SU(1,n-1)$. The reason is that the canonical symplectic form 
$\o$ on $\bR^{2n}$ is nondegenerate on the real subspace $\Pi'$ of 
$\bR^{2n}$ which corresponds to the complex subspace of $\bC^n$
defined by  $D^Aw_A=V^Aw_A=0$ under the isomorphism  $(a^A,b_B) \mapsto
(w_A)$. In fact, $\Pi'$  is complementary in $\Pi^{\perp,\omega}$ 
to the plane $\Pi\subset 
\Pi^{\perp,\omega}\subset \bR^{2n}$ spanned by the real and 
imaginary part of the complex vector
$(D^A,C_B)$.    The plane $\Pi$ is precisely the kernel of $\o$ on
$\Pi^{\perp,\omega}$. (Note that for the same reason $G^\wedge$ is not 
isomorphic to the Iwasawa subgroup
of $\SU(1,n)$.)
The complex submanifold $S':= xG'\subset S=xG^\wedge$ intersects  
all the orbits of the vector group $A$ generated by the two Killing vector
fields $\xi_1$ and $\xi_2$ transversally and exactly in one point, as follows 
from (\ref{e:holo2}).
Therefore, it is biholomorphic to the  quotient $A\setminus S$, which
is the fiber of the holomorphic submersion $\pi' : M'\ra M_{\rm sk}^\wedge$. 
This proves that the fiber 
is biholomorphic to $G'$ endowed with a left-invariant complex structure
$J'=J'(p)$, $p\in M_{\rm sk}^\wedge$. 
Using the fact that $G^\wedge$ and, hence, $G'\subset G^\wedge$ 
normalizes $A$ in $G$, one can show that the fiber metric 
corresponds to a left-invariant metric $g'=g'(p)$ on $G'$.  
Since $N_{AB}dw_Ad\bar{w}_B <0$ on $\{D^Aw_A=V^Aw_A\}\cong 
\bC^{n-2}\subset \bC^n$ we get that $(G',J',g')\cong \bC H^{n-1}$. 

In order to make the above argument more explicit, let us 
compute the K\"ahler metric of the fiber of $M'\ra M_{\rm sk}^\wedge$ in holomorphic 
coordinates and show that it is indeed the complex hyperbolic metric 
of constant  holomorphic sectional curvature $-4$. The metric of $M$ is given by \eqref{e:metric_cmap}. Let us recall
that $(\Pi^{b}_A)$ is the matrix which represents the
projection $TM_{\rm ask}\ra TM_{\rm sk}$ with respect to the 
special holomorphic coordinate frame on $M_{\rm ask}$ and a
unitary frame on $M_{\rm sk}$. By the definition of the projective special
K\"ahler metric,  we have
\be \label{e:Pi} \frac{N_{AB}}{ZN\bar{Z}} = -\delta_{ab}
\Pi^{a}_A\Pi^{b}_B + 
\frac{N_{AC}\bar{Z}^CN_{BD}Z^D}{(ZN\bar{Z})^2}, \ee
where $ZN\bar{Z}=\sum Z^AN_{AB}\bar{Z}^B=\frac{e^{-K}}{2}$, cf.\ \re{e:Kpot}. 
Multiplying \re{e:Pi} with the inverse matrix of the left hand side 
yields    
\be \label{e:Pi2}
\delta_A^B = -\tfrac12 e^{-K} \Pi_{A b} \bar \Pi_C^{~b} N^{-1\,CB} + 2 e^{K} N_{AC} \bar Z^C Z^B \ .
\ee 
Restricting to the fiber over a point $p \in M_{\rm sk}^\wedge$ and
using the identity \re{e:Pi2} we find 
\begin{displaymath}
\begin{aligned}
g_{\rm fib} = & \diff \phi^2 + \tfrac{1}{4} e^{4 \phi} |\diff \tilde \phi + b_A \diff a^A - a^A \diff b_A|^2 - \tfrac{1}{2} e^{2\phi} (\diff b_A - \bar F_{AC} \diff a^C) N^{-1\, AB} (\diff b_B - F_{BD} \diff a^D)\\ & + 2 e^{K+2\phi} |Z^A(\diff b_A - F_{AB} \diff a^B)|^2 \ ,
\end{aligned}
\end{displaymath}
which is the canonical metric of $F(n+1)$. 
Using the coordinates 
\re{holcoordEqu} the fiber metric $g_{\rm fib} := g|_{\pi^{-1}(p)}$ 
takes the following form:
\begin{equation}\label{g:fiber}
g_{\rm fib} = \tfrac{1}{4} e^{4 \phi} |\diff w^0 -2 \iu (\Im w)_A N^{-1\,AB} \diff w_B|^2 - \tfrac{1}{2} e^{2\phi} \diff \bar w_A N^{-1\, AB} \diff w_B  + 2 e^{K+2\phi} |Z^A\diff w_A|^2 \ .
\end{equation}
The metric $h_{\rm fib}:=h|_{M_p'}$ of the fiber $M_p' := (\pi')^{-1}(p)$ of 
$\pi' : M'\ra M_{\rm sk}^\wedge$ 
is obtained by first 
restricting  $g_{\rm fib}$ to the submanifold $N_p:=\pi_N^{-1}(p)\subset 
\pi^{-1}(p)$ 
defined by $D^A w_A = \tilde C$ and then 
taking the quotient by the isometric $\bR^2$-action generated by the 
Killing vector fields $\xi_1$ and $\xi_2$.  These vector fields 
can be combined in the holomorphic vector field  
\begin{displaymath}
 k:= \xi_2+i\xi_1 =
-i (\xi_1 -i\xi_2) = \bar D^A N_{AB} \frac{\partial}{\partial w_B} - 2 \iu \bar D^A (\Im w)_A \frac{\partial}{\partial w^0} \ ,
\end{displaymath}
see \re{e:holo}--\re{e:holo2}. 
Since the quotient map $\tau : N \ra M'=N/A$ is a Riemannian submersion, 
as is its restriction $\tau_p : N_p \ra M_p'$, 
the metric $h_{\rm fib}$ on $M_p'\cong F(n-1)$ is determined by the 
degenerate symmetric tensor field  
\begin{displaymath}
\begin{aligned}
(\tau_p)^*h_{\rm fib} = & \tilde g_{\rm fib} - \frac{\tilde g_{\rm fib}(k,\cdot)\tilde  g_{\rm fib}(\bar k, \cdot)+ \tilde g_{\rm fib}(\bar k,\cdot) \tilde g_{\rm fib}(k, \cdot)}{\tilde g_{\rm fib}(k,\bar k)} \ ,
\end{aligned}
\end{displaymath}
where $\tilde g_{\rm fib}= g|_{N_p}=g_{\rm fib}|_{N_p}$. 
Since $\sum Z^AN_{AB}\bar{Z}^B=\frac{e^{-K}}{2}> 0$ and therefore 
$\sum Z^AN_{AB}\bar{D}^B\neq 0$, we see as above  that 
the equivalence classes 
$[w^0,w_A]$ corresponding to 
the holomorphic $\bC$-action generated by $k$ each
contain exactly one representative which fulfills
\begin{displaymath}
 Z^A w_A = 0 \ .
\end{displaymath}
{}Recall that the index $A$ runs from $1$ to $n$. 
In particular, the $n+1$ holomorphic fiber coordinates are 
$(w^0,w_1,\ldots ,w_n)$. 
By a linear change of special coordinates $(Z^A)$, 
if necessary, we can assume that at our base point 
$p$ we have $Z^1=1$.   
Because $(Z^A)$ always has positive 
 norm and $(D^A)$ is null, we know that 
$D^a \ne D^1 Z^a$ for some $a\in \{2,\ldots ,n\}$, let us say for $a=n$.
Therefore, we can find coordinates $\{ x^0, x_a\}, a=2, \dots, n-1$, for 
the fiber of $M'\ra M_{\rm sk}^\wedge$ as follows. We put $\alpha=1/(D^n-D^1Z^n)$
and observe that the map 
\be\label{e:affineparam}
 (x^0, x_a)\mapsto (w^0,w_A) 
= \left(x^0, \alpha ((Z^n D^a - Z^a D^n) x_a -Z^n\tilde C), 
x_a, \alpha (\tilde C - (D^a - D^1 Z^a) x_a) \right) \ ,
\ee
is an affine isomorphism from $\bC^{n-1}$ onto the 
affine subspace $E\subset \bC^{n+1}$ defined by $D^Aw_A=\tilde{C}$ and 
$Z^Aw_A=0$. Therefore, it induces 
a biholomorphic map from an open subset of $\bC^{n-1}$
onto $M_p'= N_p/A \cong E\cap N_p$. 
On the complex hypersurface $\mathcal{H}:= E\cap N_p \subset N_p$ (defined 
by $Z^A w_A = 0$) we have 
\begin{displaymath}
 Z^A \diff w_A = 0 \ , \qquad D^A \diff w_A = 0 \ .
\end{displaymath}
From this one computes 
\begin{displaymath}
g(k, \cdot )|_{\mathcal{H}} =\tilde g_{\rm fib}(k,\cdot)|_{\mathcal{H}} = 0 \ ,
\end{displaymath} 
and therefore concludes that the projection $N_p \ra M_p'$ restricts to 
a biholomorphic isometry  $\mathcal{H} \ra M_p'$. 
Using this isomorphism, the metric $h_{\rm fib}$ of $M_p'$ is identified with 
the metric $g_{\mathcal{H}}=g|_{\mathcal{H}}=g_{\rm fib}|_H$ of the hypersurface $\mathcal{H}
\subset N_p$, which is 
\begin{displaymath}
\begin{aligned}
g_{\mathcal{H}} = & \tfrac{1}{4} e^{4 \phi} |\diff x^0 -2 \iu ( \tilde N_0^{b} + (\Re x)_a \tilde N_1^{ab}+(\Im x)_a \tilde N_2^{ab}) \diff x_b|^2 - \tfrac{1}{2} e^{2\phi} \diff \bar x_a \tilde N^{ab} \diff x_b \ ,
\end{aligned}
\end{displaymath} 
where 
\begin{displaymath} \begin{aligned}
 \tilde N^{ab} = & N^{-1\, ab} + \bar \alpha (\bar Z^n \bar D^a - \bar Z^a \bar D^n) N^{-1\, 1b} + \alpha N^{-1\, a1} (Z^n D^b - Z^b D^n)\\ & + |\alpha |^2 (\bar Z^n \bar D^a - \bar Z^a \bar D^n)N^{-1\, 11}( Z^n D^b - Z^b D^n)  - \bar \alpha (\bar D^a - \bar D^1 \bar Z^a) N^{-1\, nb}\\ & - \alpha N^{-1\, an}(D^b - D^1 Z^b) + |\alpha |^2 (\bar D^a - \bar D^1 \bar Z^a) N^{-1\, nn} (D^b - D^1 Z^b)\\ & - |\alpha |^2 (\bar Z^n \bar D^a - \bar Z^a \bar D^n) N^{-1\, 1n} (D^b - D^1 Z^b) \\ & - |\alpha |^2 (\bar D^a - \bar D^1 \bar Z^a) N^{-1\, n1} (Z^n D^b - Z^b D^n) 
\end{aligned}\end{displaymath}
is Hermitian and negative definite and the other coefficients are given by
\begin{displaymath} \begin{aligned}
 \tilde N_0^{a} = & \Im( \alpha \tilde C) \big( (N^{-1\,na}-Z^nN^{-1\,1a}) + \alpha (N^{-1\,n1}-Z^nN^{-1\,11}) (Z^n D^a -Z^a D^n)\\ & + \alpha (N^{-1\,1n}-N^{-1\,nn}) (D^a - D^1 Z^a) \big)\\ 
 \tilde N_1^{ab} = &  \Im( \alpha (Z^n D^a - Z^a D^n)) (N^{-1\,1b}+ \alpha N^{-1\,11}(Z^n D^b -Z^b D^n)-  \alpha N^{-1\,1n}(D^b - D^1 Z^b ) )\\ & - \Im( \alpha (D^a - D^1 Z^a)) (N^{-1\,nb}+ \alpha N^{-1\,n1}(Z^n D^b -Z^b D^n)-  \alpha N^{-1\,nn}(D^b - D^1 Z^b ) ) \\ 
  \tilde N_2^{ab} = & \Re( \alpha (Z^n D^a - Z^a D^n)) (N^{-1\,1b}+ \alpha N^{-1\,11}(Z^n D^b -Z^b D^n)-  \alpha N^{-1\,1n}(D^b - D^1 Z^b ) )\\ & +(N^{-1\,ab}+ \alpha N^{-1\,a1}(Z^n D^b -Z^b D^n)-  \alpha N^{-1\,an}(D^b - D^1 Z^b ) ) \\ & - \Re( \alpha (D^a - D^1 Z^a)) (N^{-1\,nb}+ \alpha N^{-1\,n1}(Z^n D^b -Z^b D^n)-  \alpha N^{-1\,nn}(D^b - D^1 Z^b ) ) \ .
\end{aligned}\end{displaymath}
By a linear change of holomorphic 
coordinates we can assume that $\tilde{N}^{ab}= -\d^{ab}$. Finally, 
by changing the coordinate $x^0$ into $x^0-2i\tilde{N}^a_0x_a
-2i x_a\tilde{N}^{ab}_1x_b$   
we obtain the form 
 \begin{displaymath}
\begin{aligned}
 g_{\mathcal{H}} = & \tfrac{1}{4} e^{4 \phi} |\diff x^0 -2 \iu ( (\Im x)_a M^{ab}) \diff x_b|^2 + \tfrac{1}{2} e^{2\phi} \diff \bar x_a \delta^{ab} \diff x_b \ ,
\end{aligned}
\end{displaymath} 
where 
$M^{ab}=\tilde{N}_2^{ab}-i\tilde{N}_1^{ab}=\tilde{N}^{ab}=-\d^{ab}$. 
Note that this metric has the same form as the fiber of 
$M\ra M_{\rm sk}$, which we already know has constant holomorphic 
sectional curvature.  To compare the metrics it suffices to 
put $N_{AB}=-\eta_{AB}$ (the Minkowski scalar product) and 
$(Z^A)=(1,0,\cdots ,0)$ in \re{g:fiber},
which yields 
\begin{displaymath}
\begin{aligned}
 g_{\rm fib} = & \tfrac{1}{4} e^{4 \phi} |\diff w^0 -2 \iu ( (\Im w)_A \eta^{AB}) 
\diff w_B|^2 + \tfrac{1}{2} e^{2\phi} \diff \bar w_A \delta^{AB} \diff w_B \ .
\end{aligned}
\end{displaymath}  
Changing the coordinate $w_1$ to $\bar{w_1}$ brings this metric to the 
more standard form  \re{e:metric_fibre}, but in $n+1$ instead of 
$n-1$ complex dimensions. 
\qed 

\noindent
{\bf Remark:} \label{completenesspageref} 
The above proof shows that the quotient K\"ahler manifold
$M'$ can be described as follows. As a smooth manifold, 
\[ M'=M^\wedge_{\rm sk}\times G',\] 
where $G'$ is the Iwasawa subgroup of $\SU (1,n-1)$. 
The  K\"ahler structure $(J_{M'},g_{M'})$ of $M'$ is of the form
\[ J_{M'} = J_{M^\wedge_{\rm sk}} + J',\quad g_{M'}= g_{M^\wedge_{\rm sk}} +g',\]
where $(J'(p),g'(p))_{p\in M^\wedge_{\rm sk}}$ is a family of 
left-invariant K\"ahler structures
on $G'$ such that $(G',J'(p),g'(p))$ is isomorphic to $\bC H^{n-1}$
with its standard K\"ahler structure for all $p$. 
Applying Theorem 2 of \cite{CMX}, this shows, in particular, that 
$M'$ is complete if the submanifold $M_{\rm sk}^\wedge\subset M_{\rm sk}$ 
is complete. 

We shall now consider some explicit examples of the new quotient construction applied to quaternionic K\"ahler manifolds in the image of the c-map.   

\subsubsection{Quadratic prepotential}\label{QuadprepotSec}
Let us first analyze the case of a quadratic prepotential $F$, i.e.\
$F(Z^1,\ldots ,Z^n)$ is a quadratic polynomial such that the 
real symmetric matrix $N_{AB}= \Im F_{AB}$ is of signature $(1,n-1)$. 
The corresponding $4n$-dimensional quaternionic K\"ahler manifold is 
the Hermitian symmetric space 
\begin{displaymath}
M = \frac{\U(2,n)}{\U(2)\times \U(n)} \ .
\end{displaymath}
\bp \label{Propbihquadr} In the case of quadratic prepotential, the holomorphic submersion 
$\pi : M \ra M_{\rm sk}= H_\bC^{n-1}$ of Proposition \ref{biholProp}
is a trivial holomorphic fiber bundle and $(M,J_3)$ is
biholomorphic to $H_\bC^{n-1}\times F(n+1)$. 
\ep 

\pf Since $F_{AB}$ is constant,  the fiber preserving open  
embedding $\Psi :  M \ra M_{\rm sk}\times \bC^{n+1}$ defined in \re{holcoordEqu} 
is a biholomorphic isomorphism onto its image  $M_{\rm sk}\times F(n+1)$. 
\qed 

In this case, the first condition in \eqref{e:definition_N} is automatically 
satisfied at every point of $M$ as soon as it is satisfied at one point. Hence, 
$N$ is of dimension $4n-2$ and $M'$ is of dimension $4n-4$.
\bp In the case of quadratic prepotential, 
the holomorphic submersion $\pi_N : N \ra M_{\rm sk}^\wedge=M_{\rm sk}$ of
Proposition \ref{NbiholProp} is a trivial holomorphic fiber 
bundle and the complex submanifold 
$N\subset (M,J_3)$ is biholomorphic to 
$M_{\rm sk}\times \bC \times F'(n-1) = H_\bC^{n-1}\times 
\bC\times F'(n-1)$, for any choice of null vector $(D^A)\in \bC^{1,n-1}$ 
and any $\tilde C\in \bC$. 
\ep 

\pf Since $N_{AB}$ is now constant, it follows immediately from the proof of 
Proposition \ref{NbiholProp} that the submanifold $N\subset M\cong 
M_{\rm sk}\times \bC^{n+1}$ is biholomorphic to $M_{\rm sk}\times \bC \times 
F'(n-1)$. 
\qed 
\bt \label{quadrThm} The K\"ahler manifolds $M'$ obtained from the 
quotient construction  of Corollary~\ref{cor:main} applied to 
the quaternionic K\"ahler manifold $M = \frac{\U(2,n)}{\U(2)\times \U(n)}$   
are isomorphic to 
\[ H_\bC^{n-1}\times H_\bC^{n-1},\]
for any choice of null vector $(D^A)\in \bC^{1,n-1}$ and any $\tilde C\in \bC$.  
\et 
\pf The holomorphic submersion $M'\ra M_{\rm sk}^\wedge$ 
of Theorem \ref{ThmM'} is, in this case, 
a trivial holomorphic fiber bundle over $M_{\rm sk}^\wedge=
M_{\rm sk}= H_\bC^{n-1}$. This follows from  
the proof of Theorem \ref{ThmM'},
since the constructions are now independent of $p\in M_{\rm sk}$.   
For the same reason, the metric is the product of the metric on the 
base and the metric on the fiber. 
\qed 

\subsubsection{Cubic prepotential}\label{sect:cubic}
 
Now let us turn to the case of a cubic prepotential, i.e.\
\begin{equation}\label{cubicp}
F = \tfrac16 d_{ijk} \frac{Z^i Z^j Z^k}{Z^0} \ ,
\end{equation}
where the lower case indices run from $1$ to $n-1$. Note that, from now on, 
the special coordinates $Z^I$ run from $Z^0$ to $Z^{n-1}$. 
Putting 
$z^i=Z^i/Z^0$,  the first equation in \eqref{e:definition_N} turns into
\begin{equation}\label{e:C_cubic}
C_I = \left( \begin{aligned}  d_{ijk} (\tfrac13 D^0 z^i -\tfrac12 D^i) z^j z^k \\ - d_{ijk} (\tfrac12 D^0 z^j - D^j)  z^k \end{aligned} \right) \ ,
\end{equation}
which defines a K\"ahler submanifold $M^\wedge_{\rm sk}$ of $M_{\rm sk}$ under
our general assumptions on the rank of the matrix \re{equ:GAB}, see the remark 
on page \pageref{GAB}. By means of the coordinates $z^1,\ldots,z^{n-1}$
we will identify $M_{\rm sk}$ with an open subset of $\bC^{n-1}$.  

\bp \label{rankProp2} 
Let $z_0\in M_{\rm sk}\subset \bC^{n-1}$ be a solution of the 
equation \re{e:C_cubic} and $U\subset M_{\rm sk}$ an open neighborhood of $z_0$. 
Suppose that the rank of the matrix 
\be m_{ij} := d_{ijk}(D^k-D^0z^k)\ee
is constant on $U$. Then $M^\wedge_{\rm sk}\subset M_{\rm sk}$ is a complex 
submanifold of complex codimension $r=\mathrm{rk} (m_{ij})$. 
(More generally, it suffices to assume that the rank 
of $(m_{ij})$ is constant on a  complex submanifold containing
the algebraic subset of $U\subset \bC^{n-1}$ defined by \re{e:C_cubic}.) 
\ep 
\pf The Jacobi matrix of the map $z\mapsto D^JF_{IJ}|_{Z=(1,z)}$
is given by 
\be \left( \begin{array}{c} -m_{jk}z^k\\
m_{ij}\end{array}\right).\ee
Since the first row is a linear combination of the other rows, the rank of that
matrix coincides with the rank of $(m_{ij})$.  
\qed 

\noindent
{\bf Remark:} Note that, as in the case of general prepotential, 
given a null vector $(D^I)$ at $Z=(1,z)\in  M_{\rm ask}$ we can define 
$(C_I)$ such that \re{e:C_cubic} holds at $z$. 
Therefore, we 
can always assume that $M^\wedge_{\rm sk}\neq \emptyset$.  For generic 
$(d_{ijk})$, $(D^I)$ and $z$ the rank
of $m_{ij}$ is maximal and so $\dim M^\wedge_{\rm sk}=0$.
Let us also keep in mind the trivial fact that for any $z_0=(z^i_0)
\in M_{\rm sk}$ there is always a nonzero vector $(D^I)$ 
which satisfies \eqref{e:null}
at $Z_0=(1,z_0)$. The set of all such vectors (the null cone without 
its origin) is a $\bC^*$-invariant real hypersurface of $T_{Z_0}M_{\rm ask}=\bC^n$. 
Finally, let us point out that the constant vector $(D^I)$ defining
the submanifold $M^\wedge_{\rm sk}$ is a null vector not only at $Z_0$ but 
at any point $Z$ of $M^\wedge_{\rm ask}$, since 
$$ \bar D^I N_{IJ}(Z) D^J = \frac{1}{2i}( \bar D^I C_I  - D^I \bar C_I)  =
\bar D^I N_{IJ}(Z_0) D^J  = 0 \ .$$

The following proposition can be used in explicit examples to obtain an 
upper bound on the dimension of 
$M^\wedge_{\rm sk}$, which is defined by \eqref{e:C_cubic}.
\bp \label{p:cubic} Let $z_0$ be any point of $M^\wedge_{\rm sk}$. A necessary 
condition for a vector $\alpha = \alpha^i \partial_{z^i} \in T_{z_0} M_{\rm sk}$ 
to be tangent to $M^\wedge_{\rm sk}$ is to satisfy the following equations: 
\begin{equation}
\label{e:tangent_lin}
d_{ijk} (D^j - D^0 z^j_0) \alpha^k = 0 
\end{equation} 
and  
\begin{equation}
\label{e:tangent_cub}
 d_{ijk} \alpha^i \alpha^j \alpha^k = 0 \ . 
\end{equation} 
\ep

\pf
Consider a complex analytic curve $\tau \mapsto z(\tau) = (z^i (\tau))$ 
in $M^\wedge_{\rm sk}$ through $z_0=(z^i_0)$: 
$$ z^i (\tau) = z^i_0 + \tau \alpha^i + \tau^2 \beta^i + \tau^3
\gamma^i + \dots \ .$$
Then the last $n-1$ equations of 
\eqref{e:C_cubic} are satisfied up to cubic order in $\tau$ if and only if: 
\begin{eqnarray} 0 &=& d_{ijk} (D^j - D^0 z^j_0) \alpha^k = 0 \ ,\\
\label{e:1st}
 0 &=& d_{ijk} (D^j - D^0 z^j_0) \beta^k - \tfrac12 D^0 d_{ijk} \alpha^j
\alpha^k  \quad\mbox{and}\\ \label{e:2nd}
0 &=& d_{ijk} (D^j - D^0 z^j_0) \gamma^k -D^0 d_{ijk} \alpha^j \beta^k  \ .
\label{e:3rd}
\end{eqnarray}
The first equation already gives  \eqref{e:tangent_lin}. 
Considering the $\tau^3$-component of the first equation of \eqref{e:C_cubic}
we also obtain 
\be \label{e:last} 
-d_{ijk}z^i_0 (D^j - D^0 z^j_0) \gamma^k +(2 D^0 z^i_0 - D^i  )
d_{ijk} \alpha^j \beta^k + \tfrac13 D^0 d_{ijk} \alpha^i \alpha^j
\alpha^k =0  \ .\ee
Inserting \re{e:1st}--\re{e:3rd} into \re{e:last}, we 
find $$ d_{ijk} \alpha^i \alpha^j \alpha^k  = 0 ~.$$\qed

\noindent {\bf Remark:} Note that one can always find $D^I$ s.t.\ $$d_{ijk} (D^j - D^0 z^j_0)
\alpha^k = 0$$ is not fulfilled for any $\alpha\neq 0$. On the other hand,
depending on the particular form of $d_{ijk}$, one can adjust $D^I$ in order
to obtain examples for which $\dim N$ is large. We will discuss such examples 
in the remainder of this paper. 

\noindent\textbf{A low-dimensional example}

We shall now see that  a simple low-dimensional example with a one-dimensional manifold 
$M^\wedge_{\rm sk}$ is provided by the $STU$ model with two coordinates 
fixed. The corresponding quaternionic K\"ahler manifold is the 
symmetric space 
\begin{displaymath}
M = \frac{\SO_0(4,4)}{\SO(4)\times \SO(4)}\ ,
\end{displaymath} 
which is the c-map image of the special K\"ahler manifold 
\begin{displaymath} 
M_{\rm sk} = \left(\frac{\SU(1,1)}{\U(1)}\right)^3 \ .
\end{displaymath} 
Choosing appropriate inhomogeneous coordinates  
$z=(z^1=S, z^2 =T, z^3 = U)$, the prepotential \eqref{cubicp} is determined by 
\begin{displaymath}
F(z^0=1,z)=STU \ . 
\end{displaymath}

The equation \eqref{e:C_cubic} defining the submanifold $M^\wedge_{\rm sk}
\subset M_{\rm sk}$ now reads 
\begin{align}
 2 D^0 S T U - D^S T U - D^T S U - D^U S T  &= C_0 \label{e:STU0}\ , \\
 (D^0 T - D^T) (D^0 U - D^U) &= D^T D^U - D^0 C_S \ , \label{e:STU1} \\  (D^0 S - D^S) (D^0 U - D^U) &= D^S D^U - D^0 C_T \ , \label{e:STU2}\\
   (D^0 S - D^S) (D^0 T - D^T) &= D^S D^T - D^0 C_U  \ , \label{e:STU3}
\end{align}
where $(D^S,D^T,D^U)=(D^1,D^2,D^3)$ and we are assuming that $D^0\neq 0$. 
From the last three equations we already see that two of the three coordinates, say $S$ and $T$, must be fixed to the values $\langle S \rangle := \frac{D^S}{D^0}$ and $\langle T \rangle := \frac{D^T}{D^0}$ in order to keep the 
third coordinate, here $U$, free. Note that this is not possible for arbitrary 
choices of $D^S$ and $D^T$ since the coordinates have to satisfy 
$\sum_{I,J=0}^3 N_{IJ}z^I\bar{z}^J>0$. 
 Therefore, we will assume that $(1,
\langle S \rangle,\langle T \rangle)$ 
can be extended to a vector $(1, \langle S \rangle,\langle T \rangle,\langle 
U \rangle )$ spanning a complex 
line which is positive definite with respect to the pseudo-Hermitian 
metric $(N_{IJ})$.  We will call
such vectors \emph{time-like}. 
One can check that all the above equations are solved for 
\begin{displaymath}
  C_S = D^U \langle T\rangle  \ , \qquad   C_T = D^U \langle S\rangle \ ,\qquad
   C_U= D^0 \langle S\rangle \langle T\rangle \ ,\qquad   C_0 = - D^U \langle S\rangle \langle T\rangle \ , 
\end{displaymath}
with $U$ remaining arbitrary. Therefore, the coordinate $U$ parameterises 
$M^\wedge_{\rm sk}$. It is straightforward to check that for any choice of 
$D^0\neq 0$, $D^S$ and  $D^T$ as above, the null condition 
\eqref{e:null} can be satisfied by appropriately choosing $D^U$. This 
ensures $D^0 U - D^U \ne 0$ on $M^\wedge_{\rm sk}$. The latter inequality
implies that the matrix $(m_{ij})$ of Proposition 
\ref{rankProp2} has rank two, which again proves that 
$M^\wedge_{\rm sk}\subset M_{\rm sk}$ is a one-dimensional complex submanifold.
The resulting K\"ahler manifold $M'$ has complex dimension 4.

We can also consider the quantum STU model, where the prepotential is given by
\begin{displaymath}
F(1,z)=STU + \tfrac{1}{3}T^3 \ ,
\end{displaymath}
and the corresponding $6$-dimensional special K\"ahler manifold $M_{\rm sk}$ 
admits a $4$-dimensional group of automorphisms, which acts freely on
$M_{\rm sk}$, as follows from \cite{CMX}, Example 3 in Section 4.2.

Again we can try to fix the values of one or two of the variables and use the 
remaining ones as parameters. In this case, Proposition \ref{p:cubic} immediately implies that  $T$ cannot belong to the remaining parameters. Comparing with the equations \eqref{e:STU0}--\eqref{e:STU3} for the  STU model, 
we see that only the conditions \eqref{e:STU0} and \eqref{e:STU2} are modified by the extra term in the prepotential. The new version of \eqref{e:STU2} reads
$$ (D^0 S - D^S) (D^0 U - D^U) + (D^0 T - D^T)^2 = D^S D^U + (D^T)^2 - D^0 C_T \ ,$$
which together with \eqref{e:STU1} and \eqref{e:STU3} implies that 
$T$ must be fixed to some value $\langle T\rangle$. If 
$\langle T \rangle = \frac{D^T}{D^0}$, then also $S=\langle S \rangle 
= \frac{D^S}{D^0}$ and we come back to the solution for the STU model, now with
\begin{displaymath}
  C_S = D^U \langle T\rangle  \ , \qquad   C_T = D^U \langle S\rangle + 
D^T \langle T\rangle\ ,
\end{displaymath}
\begin{displaymath}
   C_U= D^0 \langle S\rangle \langle T\rangle \ ,\qquad   C_0 = - D^U \langle S\rangle \langle T\rangle -\tfrac13 D^T \langle T\rangle^2  \ . 
\end{displaymath}
Again, $U$ parameterises $M^\wedge_{\rm sk}$ and we find again a $4$-dimensional 
K\"ahler manifold $M'$ as in the STU model. If the constant 
$\langle T \rangle$ is chosen to be real, then the term $\frac{1}{3}T^3$ in the prepotential will not contribute to 
the metric of $M'$ and so we get the same K\"ahler metric as for the 
unperturbed STU model. Otherwise, the metric will change by a conformal factor
of the form $\frac{e^{-2K_0}}{e^{-2K}}=(\frac{e^{-K_0}}{e^{-K_0} + c})^2$, where 
$c=\frac{8}{3}(\Im \langle T \rangle)^3$ and $K_0$ is the K\"ahler potential
of the unperturbed STU model.

\noindent\textbf{High-dimensional examples}

We can construct examples $M^\wedge_{\rm sk}$ with high dimension 
by extending the example above to the manifold 
\begin{equation} \label{e:SO4n}
M = \frac{\SO_0(4,n)}{\SO(4)\times \SO(n)} \ ,\quad n\ge 4\ ,
\end{equation}
which is the c-map image of 
\begin{displaymath}
M_{\rm sk} = ST[2,n-2] := 
\frac{\SU(1,1)}{\U(1)} \times \frac{\SO_0(2,n-2)}{\SO(2)\times \SO(n-2)}\ . 
\end{displaymath}
The latter has complex dimension $n-1$. By appropriately choosing 
inhomogeneous coordinates the prepotential becomes 
\begin{equation}\label{e:stuy}
F(1,z)=STU + S y^\ell y^m \delta_{\ell m} \ ,
\end{equation}
where now $z=(S,T,U,y)$, $y=(y^\ell)$ and $\ell,m = 1,\ldots ,n-4$.
For this prepotential we find from \eqref{e:C_cubic} 
\begin{align}
 \nonumber
 D^0 C_0 =&  D^0 (D^0 S - D^S) TU + D^0 S (D^0 T - D^T) U - D^0 D^U S T 
-S D^\ell\delta_{\ell m} D^m\\
 \label{e:sym:cubic}
  & + D^0 (D^0 S - D^S) y^\ell \delta_{\ell m} y^m
+  S (D^0 y^\ell-D^\ell)\delta_{\ell m} (D^0 y^m-D^m) \ , \\
 \nonumber
  D^0 C_S =& D^\ell \delta_{\ell m} D^m +D^T D^U - (D^0 T - D^T) (D^0 U - D^U) \\ \label{e:sym:S}
& - (D^0 y^\ell - D^\ell)\delta_{\ell m} (D^0 y^m - D^m)\ , \\
\label{e:sym:T}
  D^0 C_T =& -(D^0 S - D^S) (D^0 U - D^U) + D^S D^U   \ ,\\
\label{e:sym:U}
 D^0 C_U  =& -(D^0 S - D^S) (D^0 T - D^T) + D^S D^T \ ,\\
\label{e:sym:m}
  D^0 C_\ell =& - 2( D^0 S - D^S)\delta_{\ell m} (D^0 y^m - D^m)+ 2 D^S \delta_{\ell m} D^m \ .
\end{align}
From the $n=4$ example we expect that at least two directions should be 
fixed to a constant. Indeed, we have to at least fix $S$ to the value 
$\langle S \rangle = \frac{D^S}{D^0}$ in order to solve the equations 
\eqref{e:sym:T}-\eqref{e:sym:m}. From \eqref{e:sym:S} we then get an additional quadratic equation in the remaining coordinates that can be solved by choosing 
$$ T =\frac{D^T}{D^0} - \frac{\delta_{\ell m} (D^0 y^\ell -2D^\ell) y^m}{D^0 U - D^U}$$
and $(D^U,D^\ell )$ as usual such that \eqref{e:null} holds
at some base point $Z_0=(1,z_0)$ and $D^0 U - D^U \ne 0$ on $M_{\rm sk}$
(one may have to replace $M_{\rm sk}$ by a neighborhood of $z_0$ 
for the latter). The solution for $C_I$ is given by
\begin{displaymath}\begin{aligned}
  C_S =  \frac{D^T D^U}{D^0}  \ , \qquad   C_T = D^U \langle S\rangle\ ,\qquad
   C_U= \langle S\rangle D^T \ ,\\  C_\ell= 2\langle S \rangle \delta_{\ell m}D^m \ , \qquad C_0 = -  \frac{D^T D^U}{D^0} \langle S\rangle  \ . 
\end{aligned} \end{displaymath}
Hence, the dimension of $M^\wedge_{\rm sk}$ is $2(n-3)$. The manifold 
$M'$ therefore has dimension $4(n-2)$, which is eight smaller than the dimension of $M$. 
Note that the dimension of the submanifold 
$M^\wedge_{\rm sk} \subset M_{\rm sk}$ is only so high because  
we are fixing the direction S, which appears in both parts of the 
direct product manifold $M_{\rm sk}$, i.e.\ in both terms in \eqref{e:stuy}. 
This is already suggested by Proposition~\ref{p:cubic}, which implies 
that in each monomial of $\sum_{i=1}^{n-1}d_{ijk}z^iz^jz^k$ at least one 
variable must be fixed. Also, it is known that the only special 
K\"ahler manifolds which are
decomposable as a product are the symmetric spaces 
$ST[2,\ell ]$, $\ell\ge 1$, \cite{FVP}. The next step is to 
study special K\"ahler manifolds that are not symmetric and,
hence, are not decomposable. 

\noindent\textbf{Examples of homogeneous manifolds}

Let us now discuss the case of a homogeneous quaternionic manifold of negative
scalar curvature that is not necessarily symmetric. 
These manifolds have been classified 
(under certain assumptions) in 
\cite{Alekseevskii:1975,Cecotti:1989,deWit:1991nm,C}. One simple class is 
the one that is in the image of the $c \circ r$ map\footnote{The r-map
is a construction of special K\"ahler manifolds, which 
was introduced by de Wit and Van Proeyen in \cite{deWit:1991nm}.  
See \cite{CMX} for a recent discussion of some of its mathematical properties.} of the hyperbolic spaces 
$$ H_\bR^{n-2}=\frac{\SO_0(n-2,1)}{\SO(n-2)} \ ,\quad n\ge 3\ , $$ 
which is defined by the 
holomorphic prepotential
$$ F(1,z) = S (S T - x^\ell \delta_{\ell m} x^m ) \ ,$$
where $z=(S,T,x)$, $x=(x^\ell )$ and the indices $\ell, m$ 
run from $1$ to $n-3$. Thus, the corresponding special
K\"ahler manifold $M_{\rm sk}$ is still of complex dimension $n-1$. 
It is known that the corresponding quaternionic K\"ahler manifold 
$M$ can be presented as a solvable Lie group
$\mathcal{T}(p)$, $p=n-3$, of rank 3 with a left invariant quaternionic 
K\"ahler structure \cite{C}.  The only symmetric space in this series
is $\mathcal{T}(0)=\frac{\SO_0(3,4)}{\SO(3)\times \SO(4)}$. We will 
consider the case $p\ge 1$. 

Inserting this prepotential into \eqref{e:C_cubic} gives
\begin{align}
\label{e:q=-1:cubic}
  C_0 =& 2(D^0 S - D^S) S T - D^T S^2 - S ( D^0 x^\ell - 2 D^\ell)\delta_{\ell m} x^m - (D^0 S - D^S) x^\ell \delta_{\ell m} x^m \ , \\
\label{e:q=-1:S}
  C_S =& -(D^0 S - 2 D^S) T - ( D^0 T - 2 D^T) S + 
( D^0 x^\ell - 2 D^\ell)\delta_{\ell m} x^m\ , \\
\label{e:q=-1:T}
  C_T =& -( D^0S -2 D^S) S \ ,\\
\label{e:q=-1:m}
  C_\ell =& (D^0 S - 2 D^S)\delta_{\ell m} x^m +\delta_{\ell m}(D^0 x^m - 2 D^m)S \ .
\end{align}
From \eqref{e:q=-1:T} we see that $S$ is always fixed, i.e.\ locally constant
on $M^\wedge_{\rm sk}$. If $S$ is fixed to some value $\langle S\rangle $ such 
that $D^0 \langle S\rangle  - D^S \ne 0$ one can conclude from \eqref{e:q=-1:m} and \eqref{e:q=-1:S} that $T$ and $x^m$ are also fixed. If $D^0 \langle S\rangle  - D^S = 0$, we find that \eqref{e:q=-1:m} does not fix any further 
coordinates but only determines the value of $C_\ell$. 
In contrast, \eqref{e:q=-1:S} reads
$$  D^\ell\delta_{\ell m}D^m + D^0 C_S -  2D^0 D^T \langle S\rangle  = ( D^0 x^\ell - D^\ell)\delta_{\ell m} ( D^0 x^m - D^m)\ , $$
which is a quadratic equation and fixes one of the complex degrees of freedom
(which we will simply call moduli). Therefore, the 
minimal number of fixed moduli is two in the base space and two in the fiber. 
Now let us turn to the cubic equation \eqref{e:q=-1:cubic}. By using $D^0 \langle S\rangle  - D^S = 0$, we can write it as 
$$D^0  C_0 + D^0 D^T \langle S \rangle^2 +  \langle S\rangle  D^\ell\delta_{\ell m}D^m =  \langle S\rangle  ( D^0 x^\ell - D^\ell)\delta_{\ell m} ( D^0 x^m - D^m)  \ , $$ which reduces to the above quadratic equation if $C_0$ is chosen properly.
This shows that we can construct examples such that the final 
K\"ahler manifold $M'$ has complex dimension $2n-4$.

Now let us turn to a second series of homogeneous 
quaternionic K\"ahler manifolds $\mathcal{W}(p,q)$, 
which is a generalization of \eqref{e:SO4n} 
and has the prepotential
$$ F(1,z)=F(1,S,T,U,x,y)=STU + S y^\ell \delta_{\ell m} y^m+ T x^a 
\delta_{ab} x^b \ ,$$
where $x=(x^\ell)\in \bR^p$, $y=(y^a)\in \bR^q$. The 
Alekseevsky spaces $\mathcal{W}(p,q)$ are of dimension $4n=4(p+q+4)$ and are 
symmetric only if $p=0$ or $q=0$.  We will consider the case $p,q\ge 1$. 
The equation \eqref{e:C_cubic} now reads
\begin{align}
 \nonumber
 D^0 C_0 =&  D^0 (D^0 S - D^S) TU + D^0 S (D^0 T - D^T) U - D^0 D^U S T 
-S D^\ell\delta_{\ell m} D^m \\
 \nonumber
  & + D^0 (D^0 S - D^S) y^\ell \delta_{\ell m} y^m+  S (D^0 y^\ell-D^\ell)\delta_{\ell m} (D^0 y^m-D^m) \\
 \label{e:q=0:cubic}
  & +D^0 (D^0 T - D^T) x^a \delta_{ab} x^b+ T (D^0 x^a-D^a)\delta_{ab} (D^0 x^b-D^b) - T D^a\delta_{ab} D^b \ , \\
\nonumber
  D^0 C_S =& D^\ell \delta_{\ell m} D^m +D^T D^U - (D^0 T - D^T) (D^0 U - D^U) \\\label{e:q=0:S}
   & - (D^0 y^\ell - D^\ell)\delta_{\ell m} (D^0 y^m - D^m)\ , \\
\nonumber
  D^0 C_T =& D^a \delta_{ab} D^b + D^S D^U  - (D^0 S - D^S) (D^0 U - D^U) \\ \label{e:q=0:T}
  & - (D^0 x^a - D^a)\delta_{ab} (D^0 x^b - D^b) \ ,\\
\label{e:q=0:U}
 D^0 C_U  =& -(D^0 S - D^S) (D^0 T - D^T) + D^S D^T \ ,\\
\label{e:q=0:m}
  D^0 C_\ell =& - 2( D^0 S - D^S)\delta_{\ell m} (D^0 y^m - D^m)+ 2 D^S \delta_{\ell m} D^m \ , \\
\label{e:q=0:a}
  D^0 C_a =& - 2( D^0 T - D^T)\delta_{ab} (D^0 x^b - D^b) +2 D^T \delta_{ab} D^b \ .
\end{align}
We see from \eqref{e:q=0:m} that the $y^m$ can only be free if $S$ is fixed to the value $\langle S\rangle  = \frac{D^S}{D^0}$. However, from \eqref{e:q=0:T} we see that then the $x^a$ must fulfill a quadratic equation. Similarly, if one does not want to fix all the moduli $x^a$, one must fix $\langle T\rangle  = \frac{D^T}{D^0}$, cf.\ \eqref{e:q=0:a}, and \eqref{e:q=0:S} gives one quadratic equation for the $y^\ell$.\footnote{Note that the alternative of fixing e.g.\ the $y^\ell$ to $ y^\ell = \frac{D^\ell}{D^0}$ reduces the set of equations to those for the case \eqref{e:SO4n}, with the same set of solutions. In that case only the fiber dimension differs from the $M'$ obtained for \eqref{e:SO4n}.} 
Let us now turn to the cubic equation \eqref{e:q=0:cubic}. We see that for $\langle S\rangle  = \frac{D^S}{D^0}$ and $\langle T\rangle  = \frac{D^T}{D^0}$, this equation reduces to the quadratic equations encountered before, giving no new constraint on the remaining moduli. 
Therefore, four moduli in the base space are fixed, which together with the two fiber directions make six fixed moduli. The dimension of the resulting K\"ahler manifold $M'$ is thus eight smaller than that of the quaternionic manifold $M$. 

\noindent\textbf{General homogeneous manifolds}

Finally, let us discuss the case of a general homogeneous space with cubic prepotential.\footnote{The case of a quadratic prepotential has been discussed above.} The prepotential for the general cubic case is given by \cite{deWit:1991nm}
$$ F = h^1 [(h^2)^2 - h^\mu \delta_{\mu\nu} h^\nu] - h^2 h^\ell \delta_{\ell m} h^m + h^\mu \gamma_{\mu \ell m} h^\ell h^m \ .$$
Here, the index $\mu$ labels $q+1$ fields while $\ell$ labels $r$ fields that 
form representations of the $(q+1)$-dimensional Clifford algebra. Accordingly, the matrices $\gamma_{\mu}$ fulfill the Clifford algebra. The special K\"ahler base of $M$ is therefore parameterised $3+q+r$ complex scalars. Thus the dimension of $M$ is $4(4+q+r)$.

The analysis of possible dimensions of the K\"ahler quotient $M'$ is done analogously to the examples discussed above. Inserting the above prepotential into \eqref{e:C_cubic}, one finds
\begin{align}
 \nonumber
 D^0 C_0 =& 2 D^0 h^1 h^2(D^0 h^2 - D^2) - D^0 D^1  ( (h^2)^2 - h^\mu \delta_{\mu\nu} h^\nu) - 2 D^0 h^1 h^\mu \delta_{\mu\nu}(D^0 h^\nu - D^\nu) \\
 \nonumber
  & -  D^0 (D^0 h^2 - D^2) h^\ell \delta_{\ell m} h^m - h^2  (D^0 h^\ell-D^\ell)\delta_{\ell m} (D^0 h^m-D^m) \\
 \nonumber
  &  + D^0 \gamma_{\mu \ell m} (D^0 h^\mu - D^\mu)h^\ell h^m + \gamma_{\mu \ell m} h^\mu (D^0 h^\ell-D^\ell) (D^0 h^m-D^m) \\  \label{e:genhom:cubic}
   & + h^2 D^\ell\delta_{\ell m} D^m - h^\mu \gamma_{\mu \ell m} D^\ell D^m \ , \\
\label{e:genhom:1}
  D^0 C_1 =& (D^2)^2  - D^\mu \delta_{\mu\nu} D^\nu -(D^0 h^2 - D^2)^2 + (D^0 h^\mu - D^\mu)\delta_{\mu\nu} (D^0 h^\nu - D^\nu)\ , \\
\nonumber
  D^0 C_2 =& 2 D^1 D^2 - D^\ell \delta_{\ell m} D^m - 2(D^0 h^1 - D^1) (D^0 h^2 - D^2) \\ 
  \label{e:genhom:2} & + (D^0 x^\ell - D^\ell)\delta_{\ell m} (D^0 x^m - D^m) \ ,\\
\nonumber
 D^0 C_\mu  =&  \gamma_{\mu \ell m}D^\ell2 D^m - 2 D^1 \delta_{\mu \nu} D^\nu + 2 (D^0 h^1 - D^1) \delta_{\mu\nu}(D^0 h^\nu - D^\nu) \\    
  \label{e:genhom:mu} & -  (D^0 h^\ell - D^\ell) \gamma_{\mu \ell m} (D^0 h^m - D^m) \ ,\\
\nonumber
  D^0 C_\ell =& 2[( D^0 h^2 - D^2)\delta_{\ell m} - ( D^0 h^\mu - D^\mu)\gamma_{\mu \ell m} ] (D^0 h^m - D^m)\\ \label{e:genhom:m} & - 2 D^2 \delta_{\ell m} D^m+ 2 D^\mu \gamma_{\mu \ell m} D^m \ .
\end{align}
From \eqref{e:genhom:m} we see that the only $h^\ell$ that can stay massless are those in the kernel of the matrix 
$$M_{\ell m}(h^2, h^\mu) = [( D^0 h^2 - D^2)\delta_{\ell m} + ( D^0 h^\mu - D^\mu)\gamma_{\mu \ell m} ] \ .$$
On the other hand, a direction in the $(h^2,h^\mu)$-plane can only remain unfixed if $D^0 h^\ell - D^\ell=0$ holds for at least some of the scalars $h^\ell$. In general, the minimal set of fixed scalars consists of just $(h^2, h^\mu)$. If we fix these scalars to $\langle h^2\rangle  = \frac{D^2}{D^0}$ and $\langle h^\mu\rangle  = \frac{D^\mu}{D^0}$, then \eqref{e:genhom:1} and \eqref{e:genhom:m} are fulfilled for all values of the $h^\ell$. Furthermore, we find $q+1$ quadratic equations for the $h^\ell$ from \eqref{e:genhom:2} and \eqref{e:genhom:mu}, which also solve \eqref{e:genhom:cubic}. In total, this gives $2q+2$ fixed (complex) directions, leading to a K\"ahler quotient $M'$ of (complex) dimension $2r+4$.\footnote{Alternatively, one could choose to fix $h^1$ and all $h^\ell$, with one additional constraint coming from \eqref{e:genhom:1}, resulting in a K\"ahler quotient $M'$ of complex dimension $4+2q+r$. Depending on $q$ and $r$, this might be larger or smaller than $2r+4$. Since $r$ must be a multiple of the dimension of the fundamental representation of the $q$-dimensional Clifford algebra, generically $r$ will be much larger than $q$. Note that there is also the possibility of fixing some $h^\ell$ and some $(h^2,h^\mu)$, which we do not discuss any further here.}

\section{K\"ahler quotients and spontaneous partial supersymmetry breaking} \label{sect:physics}

The construction of the K\"ahler quotient of quaternionic-K\"ahler 
manifolds presented in Section \ref{sect:quotient} first arose 
in the physics literature in the derivation of the low-energy 
effective action of 
spontaneous $\cN=2$ to $\cN=1$ supersymmetry breaking in supergravity 
\cite{JL,LST1,LST2}. 
Let us close this paper by linking the mathematical
analysis of the previous sections
to the physical perspective
of refs.~\cite{JL,LST1,LST2}. 

The spectrum of  $\cN=2$ supergravity includes the 
gravitational
multiplet together with $n_{\rm  v}$ vector- and $n_{\rm  h}$ hypermultiplets.
Each hypermultiplet contains four real scalars
which together span a  $4n_{\rm  h}$-dimensional 
field space $M$ that is constrained by
 $\cN=2$ supersymmetry to be quaternionic-K\"ahler.
A necessary condition for a maximally-symmetric solution of the
$\cN=2$ supergravity field equations to preserve only $\cN=1$
supersymmetry is that two isometries 
of the quaternionic-K\"ahler
manifold are gauged \cite{Ferrara:1995gu,LST1}. The Higgs mechanism then makes
the corresponding two vector fields massive, with the charged scalars
providing the longitudinal degrees of freedom.
Consistency with $\cN=1$ supersymmetry demands that these isometries
satisfy the assumed properties of Theorem~\ref{thm:main}.

In order to derive an effective action valid below the scale
of supersymmetry breaking $m_{3/2}$ one needs to integrate out all
fields with masses of order $m_{3/2}$.
Integrating out massive scalar fields corresponds to taking a
submanifold $N\subset M$, while integrating out the two massive 
vector fields corresponds to taking the quotient with respect 
to the two-dimensional Abelian Lie group $A$ generated by the two
Killing vectors, as specified in Theorem~1.
The two charged scalars act as  Goldstone bosons and 
are removed from the scalar field space by the quotient construction described in Theorem~\ref{thm:main}. Consistency with $\cN=1$ supersymmetry implies that the resulting scalar field space  $M'=N/A$ should be K\"ahler.

For generic  quaternionic-K\"ahler manifolds $M$
the precise identification of massive versus massless fields
or, in other words the identification of the submanifold $N$,
is difficult. However, for the case of special quaternionic-K\"ahler manifolds, i.e.
manifolds in the image of the c-map, $N$ is determined
by \eqref{e:definition_N}, which we repeat here for convenience:
\begin{equation}\label{e:lighthhm} 
D^A {F}_{AB}(Z) = C_B \ , \qquad
D^A (b_A - {F}_{AB} a^B ) = \tilde C \ .
\end{equation}
These equations give $2r_{\mathrm{F}} +2$ real conditions, where $r_{\mathrm{F}} =  {\rm rank} ({F}_{ABC}D^C )$. From this one can read off the dimension of the submanifold $N$ to be $4n_{\rm h} -2 (r_{\mathrm{F}}+1)$. 

The dimension of the quotient ${M^\prime}$ is two less than that of the submanifold $N_{\rm  h}$. Therefore, the specific dimensions of the quotient ${M^\prime}_{\rm  h}$ is model-dependent and depends on the number of hypermultiplet scalars which remain massless, i.e.\ on the dimension of $N$. The maximal rank of ${F}_{ABC}D^C $ is $n_{\rm h} - 1$ due to ${F}_{ABC}X^A = 0$, therefore for generic $F$ and $D^A$ the dimension of $N$ is $2n_{\rm h}$, cf.\ Proposition \ref{rankProp}. In other words, generically all moduli in the special K\"ahler base of the special quaternionic-K\"ahler manifold are fixed. However, only two of the axionic scalars in the $G$-fiber are fixed. For special choices of the prepotential $F$ and fine-tuned $D^A$ one can increase the dimension of ${M^\prime}$, as discussed in detail in Section \ref{sect:ex}.

$\cN=2$ gauged supergravities in four dimensions appear in the low-energy limit of compactifications of string theory on Calabi-Yau and, more generally, $SU(3)\times SU(3)$-structure manifolds. In all these theories the quaternionic-K\"ahler manifold are of the special form described in \cite{Cecotti:1988qn,Ferrara:1989ik}. In the limit of large volume the holomorphic prepotential simplifies to become cubic. In Section \ref{sect:cubic} we analysed a large class of special quaternionic-K\"ahler manifolds with cubic prepotentials, including the examples of general homogeneous manifolds classified in \cite{deWit:1991nm} and the inhomogeneous quantum STU model. We found that it is possible to obtain both high- and low-dimensional moduli spaces, with the latter being generic. From the perspective of string theory compactifications, the fact that we generically find low-dimensional moduli spaces is particularly attractive as it suggests that moduli stabilisation can be easily implemented.

\section*{Acknowledgments}
We would like to thank Antoine Van Proeyen for useful conversations.
This work was partly supported by the German Science Foundation (DFG) under the
Collaborative Research Center (SFB) 676 “Particles, Strings and the Early Universe”.
The work of P.S. was partly supported by the Swiss National Science Foundation.
The work of H.T. is supported by the DSM CEA/Saclay, the ANR grant 08-JCJC-0001-0
and the ERC Starting Independent Researcher Grant 240210 - String-QCD-BH.

\end{document}